\documentclass[11pt]{amsart}
\usepackage{amssymb, amsmath, amscd, eucal, rotating}

\input xy
 \xyoption{all}

\numberwithin{equation}{section}

\theoremstyle{plain}

\newtheorem{proposition}{Proposition}[section]
\newtheorem{theorem}[proposition]{Theorem}		
\newtheorem{corollary}[proposition]{Corollary}
\newtheorem{lemma}[proposition]{Lemma}

\theoremstyle{definition}

\newtheorem{remark}[proposition]{Remark}

\newcommand{\CBbb}{\mathbb C}

\newcommand{\ZBbb}{\mathbb Z}
\newcommand{\QBbb}{\mathbb Q}

\newcommand{\BBbb}{\mathbb B}

\newcommand{\Gr}{\mathop{\rm Gr}\nolimits}

\newcommand{\YMH}{\mathop{\rm YMH}\nolimits}

\newcommand{\End}{\mathop{\rm End}\nolimits}

\newcommand{\Hom}{\mathop{\rm Hom}\nolimits}

\newcommand{\ad}{\mathop{\rm ad}\nolimits}

\newcommand{\dbar}{\bar\partial}
\newcommand{\lra}{\longrightarrow}

\newcommand{\pr}{\mathop{\rm pr}\nolimits}

\newcommand{\G}{\mathcal G}
\newcommand{\B}{\mathcal B}

\newcommand{\A}{\mathcal A}

\newcommand{\Tor}{{\mathcal I}}

\newcommand{\Bcal}{\mathcal{B}}

\newcommand{\Acal}{\mathcal{A}}

\newcommand{\Scal}{{\mathcal S}}
\newcommand{\Gcal}{{\mathcal G}}
\newcommand{\Mcal}{{\mathcal M}}
\newcommand{\Ncal}{{\mathcal N}}
\newcommand{\Ccal}{{\mathcal C}}

\newcommand{\SL}{\mathsf{SL}}

\newcommand{\SU}{\mathsf{SU}}

\newcommand{\PU}{\mathsf{PU}}
\newcommand{\U}{\mathsf{U}}

\newcommand{\doubleslash}{\bigr/ \negthinspace\negthinspace \bigr/}

\DeclareMathOperator{\im}{im}
\DeclareMathOperator{\rank}{rank}

\begin{document}


\title[Cohomology of $\U(2,1)$ representation varieties]
{Cohomology of $\U(2,1)$ representation \\ varieties of surface groups}

\author[Richard A. Wentworth]{Richard A. Wentworth}

\address{Department of Mathematics,
   University of Maryland,
   College Park, MD 20742, USA}
\email{raw@umd.edu}
\thanks{R.W. supported in part by NSF grant DMS-1037094}

\author[Graeme Wilkin]{Graeme Wilkin}

\address{Department of Mathematics  \\
National University of Singapore \\
Singapore 119076 }
\email{graeme@nus.edu.sg}

\thanks{}


\date{\today}

\begin{abstract} 
In this paper we use the Morse theory of the Yang-Mills-Higgs functional on the singular space of Higgs bundles 
on Riemann surfaces to compute the equivariant cohomology of the space of semistable 
$\U(2,1)$ and $\SU(2, 1)$ Higgs bundles with fixed Toledo invariant. In the non-coprime case this gives new results about the topology of the $\U(2,1)$ and $\SU(2,1)$ character varieties of  surface groups. The main results are a calculation of the equivariant Poincar\'e polynomials, a Kirwan surjectivity theorem in the non-fixed determinant case, and a description of the action of the Torelli group on the equivariant cohomology of the  character variety. This builds on earlier work for stable pairs and rank $2$ Higgs bundles.
\end{abstract}



\maketitle

\thispagestyle{empty}


\baselineskip=16pt
\setcounter{footnote}{0}


\section{Introduction}

Let $X$ be a closed Riemann surface of genus $g\geq 2$.
Choose complex hermitian vector  bundles $E_1$, $E_2$ on $X$ with $\rank E_i=i$ and 
$\deg E_i=d_i$.  Let $\Bcal(d_1,d_2)$ denote the space of $\U(2,1)$-Higgs bundle 
structures on $E_2\oplus E_1$, and let $\Gcal$ denote the 
group of $\U(2)\times\U(1)$ gauge transformations.  
For a holomorphic line bundle $\Lambda\to X$ of degree $d_1+d_2$, 
let $\Bcal_{\Lambda}(d_1,d_2)$ be the subspace defined by restricting to holomorphic structures with fixed holomorphic isomorphism
$E_1\otimes \det E_2\cong \Lambda$, and let $\Gcal_0$ denote the group  of ${\mathsf S}(\U(2)\times\U(1))$ gauge transformations. 
Denote the corresponding moduli spaces of semistable Higgs bundles by 
\begin{align}
\begin{split} \label{eqn:moduli}
\Mcal(d_1,d_2)&=\Bcal^{ss}(d_1,d_2)\doubleslash\Gcal^\CBbb\\
\Mcal_\Lambda(d_1,d_2)&=\Bcal^{ss}_{\Lambda}(d_1,d_2)\doubleslash\Gcal_0^\CBbb
\end{split}
\end{align}
The main result of this paper is a computation of  the $\Gcal$ and $\Gcal_0$-equivariant
Betti numbers of $\Bcal^{ss}(d_1,d_2)$ and
$\Bcal^{ss}_\Lambda(d_1,d_2)$.

Tensoring by line bundles and dualizing give equivariant  isomorphisms of 
these spaces. The distinct cases are therefore enumerated 
by the mod $3$ values $d_1+d_2 \equiv 0,1$, which we will refer to as 
the \emph{non-coprime} and \emph{coprime} cases,
 respectively.  The moduli spaces are nonempty only if 
 $\tau=\tau(d_1,d_2)=\tfrac{2}{3}(2d_1-d_2)$
satisfies  $|\tau|\leq 2g-2$. 
 By duality, we will assume  without loss of generality that $\tau\geq 0$.
 For a rank 2 hermitian vector bundle $E\to X$ of degree $d$,
we also introduce the space $\Ccal(E)$ of \emph{holomorphic pairs} consisting of holomorphic structures on $E$ plus a choice of holomorphic section.  Given a real number $\sigma$, $d/2\leq \sigma\leq d$, let 
$\Ccal_\sigma(E)\subset \Ccal(E)$ denote the space of 
$\sigma$-semistable pairs in the sense of Bradlow \cite{Bradlow91,
BradlowDaskal91}.  We denote the corresponding moduli space $\Ncal_\sigma(E)=\Ccal_\sigma(E)\doubleslash\Gcal^\CBbb(E)$,
where 
$\Gcal^\CBbb(E)$ is the complexification of the group $\Gcal(E)$ of unitary gauge transformations of  $E$.
 For 
  generic  $\sigma$ (generic means semistable implies stable, which occurs at noninteger values in $(d/2, d)$),
 the Poincar\'e polynomials of $\Ncal_\sigma(E)$
were computed in \cite{Thaddeus94}.  
For general values of $\sigma$  (not necessarily generic),
 the $\Gcal(E)$-equivariant cohomology of $\Ccal_\sigma(E)$
was computed  in \cite{WentworthWilkin11}.

To state the main results, set
\begin{equation} \label{eqn:sigma}
\sigma(d_1,d_2)=
2g-2+(d_2-2d_1)/3 
\end{equation}
We also let
 $J(X)$ and $S^mX$ denote the Jacobian variety and $m$-th symmetric product of $X$, respectively.
With this background we have

\begin{theorem}[$\U(2,1)$ Higgs bundles]          \label{thm:simplified-polynomials}
Fix $(d_1,d_2)$ such that $0\leq \tau(d_1,d_2)\leq 2g-2$ and 
$d_1+d_2\equiv 0\negthinspace \mod 3$. Then the $\Gcal$-equivariant Poincar\'e polynomial is given by
\begin{align} 
P_t^\mathcal{G}(\mathcal{B}^{ss}(d_1, d_2)) &=\frac{1}{(1-t^2)}  P_t^{\Gcal(E)}(\Ccal_{\sigma(d_1,d_2)}(E))P_t(J(X)) \notag  \\
 &  + \hskip-1cm \sum_{\frac{1}{3}(d_1+d_2) < \ell \leq d_2-d_1+2g-2}  \frac{t^{2(g-1+2\ell-d_2)}}{(1-t^2)}    P_t(J(X)) P_t(S^{d_2-d_1+2g-2-\ell} X) P_t(S^{d_1-\ell+2g-2} X) 
\label{eqn:non-coprime-polynomial}
\end{align}
where $\deg E=d_2-2d_1+4g-4$.
For $d_1+d_2\equiv 1\negthinspace \mod 3$,
\begin{align}
 P_t(\mathcal{M}(d_1, d_2)) 
&= (1-t^2)P_t^\mathcal{G}(\mathcal{B}^{ss}(d_1,d_2)) 
 =
 P_t(\Ncal_{\sigma(d_1,d_2)}(E))P_t(J(X)) \notag\\
 & +\hskip-.75cm \sum_{\frac{1}{3}(d_1+d_2) < \ell \leq d_2-d_1+2g-2}  t^{2(g-1+2\ell-d_2)}   P_t(J(X)) P_t(S^{d_2-d_1+2g-2-\ell} X) P_t(S^{d_1-\ell+2g-2} X) .
\label{eqn:coprime-polynomial}
\end{align}
\end{theorem}

In order to state the result for fixed determinant, let $\widetilde S(m_1,m_2)$ denote the 
pullback by the $3^{2g}$-fold cover $J(X)\to J(X): L\mapsto L^3$ of the product $S^{m_1}X\times S^{m_2}X$, where the 
map to $J(X)$ factors through $(L_1,L_2)\mapsto L_1^\ast L_2\Lambda$.  
The Poincar\'e polynomial of $\widetilde S(m_1,m_2)$ was computed by 
Gothen \cite{Gothen02} (see also Corollary \ref{cor:gothen} below).

\begin{theorem}[$\SU(2,1)$ Higgs bundles]       \label{thm:su-simplified-polynomials}
Fix $(d_1,d_2)$ such that $0\leq \tau(d_1,d_2)\leq 2g-2$ and $d_1+d_2\equiv 0\negthinspace \mod 3$. Then the $\Gcal_0$-equivariant Poincar\'e polynomial is given by
\begin{align} 
P_t^{\Gcal_0}\left( \Bcal^{ss}_{\Lambda}(d_1,d_2) \right) & =\frac{1}{(1-t^2)}  P_t^{\Gcal(E)}(\Ccal_{\sigma(d_1,d_2)}(E))P_t(J(X)) \notag  \\
& + \hskip-1cm \sum_{\frac{1}{3}(d_1+d_2) < \ell \leq d_2-d_1+2g-2} t^{2(g-1+2\ell-d_2)}  P_t(\widetilde{S}(d_2-d_1+2g-2-\ell, d_1-\ell+2g-2) )
 \label{eqn:fixed-det-noncoprime}
\end{align}
where $\deg E=d_2-2d_1+4g-4$.
For $d_1+d_2\equiv 1 \negthinspace\mod 3$,
\begin{align}
P_t(\Mcal_\Lambda(d_1,d_2))&=
P_t^{\Gcal_0}\left( \Bcal^{ss}_{\Lambda}(d_1,d_2) \right)  = P_t(\Ncal_{\sigma(d_1,d_2)}(E)) \notag \\
& \quad + \hskip-1cm \sum_{\frac{1}{3}(d_1+d_2) < \ell \leq d_2-d_1+2g-2} t^{2(g-1+2\ell-d_2)}  P_t(\widetilde{S}(d_2-d_1+2g-2-\ell, d_1-\ell+2g-2) ) 
\label{eqn:fixed-det-coprime}
\end{align}
\end{theorem}

Eq.'s  \eqref{eqn:coprime-polynomial} and \eqref{eqn:fixed-det-coprime} have been previously 
obtained by Gothen \cite{Gothen02}. 
 In the  coprime case the moduli space is smooth, and one may use the moment map 
associated to Hitchin's $S^1$-action  as a Morse-Bott function. Critical points correspond 
to fixed points of the $S^1$-action, and the cohomology of these critical sets 
(as well as their Morse indices) can be computed. 
As outlined below, the derivation of the Poincar\'e polynomials in this paper 
is different from that of \cite{Gothen02}.
Indeed, 
 showing that the two results agree in the coprime case actually depends on the 
results of \cite{Thaddeus94, WentworthWilkin11}.
The  stable pairs moduli space that occurs in Gothen's calculations 
has a different stability parameter $\sigma$ 
 to that which occurs in the calculations of this paper, 
and one needs to look at different critical sets 
for the terms corresponding to the flips that 
relate the two different Bradlow spaces. 
Therefore the connection between the two pictures is somewhat complicated and  is not merely a comparison of critical sets.

We also point out the following special case (see Section \ref{sec:equivariant-betti-numbers}).
\begin{corollary} \label{cor:maximal}
In the maximal case $\tau(d_1,d_2)=2g-2$,
\begin{equation*}
P_t^\mathcal{G}(\Bcal^{ss}(d_1,d_2)) = \frac{1}{(1-t^2)^2} P_t(J(X))^2 
\end{equation*}
\end{corollary}

\noindent
This is exactly what one would expect from Theorem B in \cite{BGG03}.

We now describe the relationship with representation varieties.
Fix $p\in X$, and let $\pi=\pi_1(X,p)$ denote the 
fundamental group acting by deck transformations on the universal cover $\widetilde X$ of $X$.
Let $\omega_{\BBbb^2}$ denote the 
complete
$\PU(2,1)$-invariant K\"ahler metric on the complex ball
$\BBbb^2\subset \CBbb^2$,
 normalized  to have constant holomorphic sectional curvature $-1$.
Given $\rho:\pi\to \PU(2,1)$, choose a $\rho$-equivariant map
$f:\widetilde X\to \BBbb^2$.  Then $f^\ast\omega_{\BBbb^2}$ is a
$\pi$-invariant form, and the \emph{Toledo invariant} of $\rho$ is
by definition
\begin{equation} \label{eqn:toledo} 
\tau(\rho)=\frac{1}{2\pi}\int_X f^\ast\omega_{\BBbb^2}
\end{equation}
By \cite{Toledo79}, $\tau(\rho)$ is an 
integer that is constant on 
connected components of the representation variety, and  which satisfies the bound $|\tau(\rho)|\leq 2g-2$.
Extend the definition of $\tau(\rho)$ to representations of $\pi$ to 
$\SU(2,1)$ and $\U(2,1)$ by projection to $\PU(2,1)$.
Let $\Hom_\tau(\pi, G)$, $G=\SU(2,1)$, $\U(2,1)$, or $\PU(2,1)$, denote the subset of representations
$\pi\to G$ with Toledo invariant $=\tau$, and let $\Hom_\tau(\pi, G)\doubleslash G$ be the corresponding moduli space of conjugacy classes of semisimple representations.
By work of Hitchin, Simpson, Corlette and Donaldson (\cite{Hitchin87,Corlette88,Donaldson87, Simpson88}; see also \cite{BGG03}) we have
\begin{align*}
\Hom_\tau(\pi, \U(2,1))\doubleslash \U(2,1)&\simeq
\Mcal(d_1,d_2)\\
\Hom_\tau(\pi, \SU(2,1))\doubleslash \SU(2,1)&\simeq
\Mcal_\Lambda(d_1,d_2)
\end{align*}
as real algebraic varieties,
where $d_1+d_2=0$,
and the two definitions of the Toledo invariant agree:
$\tau=\tau(\rho)=\tau(d_1,d_2)$.
 As explained in \cite{DWW10}, 
the results of this paper also compute the equivariant cohomology of these representation varieties
(in this paper we take rational coefficients unless otherwise indicated). 
   
 \begin{theorem} Let $d_1+d_2=0$ and $\tau = \tfrac{2}{3}(2d_1-d_2)$.  Then there are isomorphisms of equivariant cohomologies
 \begin{align*}
 H^\ast_{\U(2,1)}(\Hom_\tau(\pi, \U(2,1)))&\simeq H^\ast_{\Gcal}(\Bcal^{ss}(d_1,d_2)) \\
  H^\ast_{\SU(2,1)}(\Hom_\tau(\pi, \SU(2,1)))&\simeq H^\ast_{\Gcal_0}(\Bcal^{ss}_\Lambda(d_1,d_2)) 
  \end{align*}
\end{theorem}

Tensoring a rank-$n$ bundle by the $n$-torsion
points in the Jacobian variety $J(M)$ leaves the determinant
unchanged.  Hence, the group $\Gamma_n=H^1(M,\ZBbb/n)$ acts on 
fixed determinant moduli spaces, and the study of its induced action
 on the  cohomology 
of moduli spaces
goes back to Harder-Narasimhan \cite{HarderNarasimhan74}.  
In terms of representations, this action corresponds to the 
different possible lifts of $\PU(n)$ bundles to $\SU(n)$. 
 More precisely,  in our situation $\Mcal_\Lambda(d_1,d_2)$ is a 
 $\Gamma_3$-covering of a connected component of
$\Hom(\pi, \PU(2,1))\doubleslash\PU(2,1)$.  Furthermore, by a theorem of Xia \cite{Xia00} 
the connected components of the space of $\PU(2,1)$ representations
are in 1-1 correspondence with the mod $3$ values of $d=\deg\Lambda$ and 
 the possible values of 
the Toledo invariant $|\tau(d_1,d_2)|\leq 2g-2$.  As in the  theorem above we have
 \begin{equation} \label{eqn:pu_cohomology}
 H^\ast_{\PU(2,1)}(\Hom_{\tau,d}(\pi, \PU(2,1)))=
\left[ H^\ast_{\Gcal_0}(\Bcal^{ss}_\Lambda(d_1,d_2))\right]^{\Gamma_3}
 \end{equation}
 where $d=d_1+d_2$, $\tau(d_1,d_2)=\tau$, and  the superscript indicates the $\Gamma_3$-invariant part
 of the cohomology.

It was shown in Atiyah-Bott \cite{AtiyahBott83}, and illustrated further in 
\cite{DWWW11} for $\SL(2,\CBbb)$, that
 the action of $\Gamma_n$ is also the key to
understanding \emph{Kirwan surjectivity}, 
which we now define.  Since the spaces $\Bcal(d_1,d_2)$ and
$\Bcal_\Lambda(d_1,d_2)$ are contractible,  the inclusions
$\Bcal^{ss}(d_1,d_2)\hookrightarrow \Bcal(d_1,d_2)$ and 
$\Bcal^{ss}_\Lambda(d_1,d_2)\hookrightarrow \Bcal_\Lambda(d_1,d_2)$
give maps
\begin{align}
\begin{split} \label{eqn:kirwan_map}
\kappa : H^\ast(B\Gcal) &\lra H^\ast_\Gcal(\Bcal^{ss}(d_1,d_2))\\
\kappa_0 : H^\ast(B\Gcal_0) &\lra H^\ast_{\Gcal_0}(\Bcal^{ss}_\Lambda(d_1,d_2))
\end{split}
\end{align}
which we call \emph{Kirwan maps}. We say that Kirwan surjectivity holds  if $\kappa$ (or $\kappa_0$) is surjective.
For 
 $\U(n)$ and $\SU(n)$  bundles, it turns out that the Kirwan maps are always surjective \cite{AtiyahBott83}.
This is a consequence of the perfection of the Harder-Narasimhan (and
Morse) stratification.
It is also the case that $\Gamma_n$ acts trivially on
$H^\ast(B\Gcal_0)$, and so surjectivity implies the same for the
cohomology of the representation varieties.
On the other hand, for $\SL(2,\CBbb)$ Higgs bundles, $\kappa_0$ is not in
general surjective (cf.\ \cite{DWWW11}).

 Continuing in this vein, we show
in this paper that a certain modification of the Harder-Narasimhan stratification 
for $\U(2,1)$ Higgs bundles is $\Gcal$-equivariantly perfect (Theorem
\ref{thm:perfection}),
and hence Kirwan surjectivity holds in this case.
We also show that $\Gamma_3$ acts trivially on the
equivariant cohomology of the moduli space of $\SU(2,1)$ Higgs
bundles if and only if Kirwan surjectivity holds. 
In the fixed determinant case,  surjectivity holds
for only
about a third of  the components.
\begin{theorem} \label{thm:kirwan}
Kirwan surjectivity holds for the moduli spaces of  
$\U(2,1)$ and $\PU(2,1)$ Higgs
bundles.  Kirwan surjectivity holds for  the moduli spaces of
$\SU(2,1)$ Higgs bundles if and only if the Toledo invariant satisfies
$|\tau|>\tfrac{4}{3}(g-1)$.
\end{theorem}

The action of $\Gamma_n$ is also closely intertwined with  the action of
the Torelli group $\Tor(X)$, defined as the subgroup of the mapping
class group that acts trivially on the homology of $X$ (see
\cite{DWW10}).  Since $\Tor(X)$ is 
a subgroup of the outer automorphism group of $\pi$, it acts on representation varieties by precomposition, and
 the induced action on  equivariant cohomology  commutes with $\Gamma_n$.  On the other hand,
by
results of Looijenga \cite{Looijenga97} characters of $\Gamma_n$ give rise to
projective unitary representations of $\Tor(M)$ over 
cyclotomic fields. 
In Theorem \ref{thm:torelli}, we explicitly determine the
representations that appear for the action of $\Gamma_3\times
\Tor(M)$ on the moduli space of $\SU(2,1)$ Higgs bundles.  As a consequence, we prove

\begin{theorem}\label{thm:Torelli}
The group $\Gamma_3\times \Tor(X)$ acts trivially on the equivariant cohomology of the moduli spaces of
$\U(2,1)$ and $\PU(2,1)$ representations of $\pi$.
The Torelli group $\Tor(X)$ $($resp.\ the group $\Gamma_3$$)$ acts
trivially on the equivariant cohomology of the moduli spaces of
$\SU(2,1)$ representations  if and only if the Toledo invariant
satisfies $|\tau|\geq \tfrac{4}{3}(g-1)$ $($resp.\ $|\tau| >
\tfrac{4}{3}(g-1)$$)$.
\end{theorem}

\noindent
The borderline case $\tau=\tfrac{4}{3}(g-1)$ (which occurs only for 
$g\equiv 1\negthinspace\mod 3$) gives further examples in higher genus of
  representation varieties where Kirwan surjectivity fails but
 where the Torelli group nevertheless acts trivially on equivariant cohomology 
(this also occurs for $\SL(2,\CBbb)$ bundles, but only  when $g=2$).
Gothen also studies the action of $\Gamma_3$ on the cohomology of the moduli space and shows in \cite[Proposition 4.2]{Gothen02} that in general it acts non-trivially on $H^*(\Mcal_\Lambda(d_1,d_2))$ in the coprime case.

The method of proof for the results above is an extension of the equivariant Morse theory techniques of Atiyah-Bott and Kirwan from \cite{AtiyahBott83} and \cite{Kirwan84} to the singular space of Higgs bundles. This continues a program begun in \cite{DWWW11, DWW10} (for rank $2$ Higgs bundles) and \cite{WentworthWilkin11} (for rank $2$ stable pairs), and we use these results as part of our calculations for the $\U(2,1)$ case.
The basic strategy is to use the Yang-Mills-Higgs functional as an equivariant Morse function on the spaces of $\U(2,1)$ Higgs bundles (resp. $\SU(2,1)$ Higgs bundles), where equivariance is defined with respect to the  group of gauge transformations in  the maximal compact subgroup of $\U(2,1)$ (resp.\ $\SU(2,1)$). 

In Section \ref{sec:stratifications} we describe the stratification of the space of Higgs bundles by the gradient flow of the Yang-Mills-Higgs functional and assert that the gradient flow of the Yang-Mills-Higgs functional on the space of $\U(2,1)$ Higgs bundles induces a Morse stratification identical to the Harder-Narasimhan stratification. Another result of this section is that the equivariant  cohomology of the critical sets can be computed inductively in terms of lower rank Higgs bundles.

The major subtlety induced by the singularities in the space of Higgs bundles occurs in the study of the change in cohomology when attaching each of the Morse/Harder-Narasimhan strata to the union of lower strata. The Morse index is not constant on each connected component of the set of critical points, and so instead of attaching a bundle over the critical set (as in the usual Morse-Kirwan theory) one has to attach a more general space that fibers over the critical set. Section \ref{sec:calculations} contains a detailed analysis of these spaces and a calculation of their cohomology.

The Poincar\'e polynomial calculations are summarized in Section \ref{sec:equivariant-betti-numbers}. The key point is that the spaces described in the previous paragraph appear as extra terms in the Poincar\'e polynomials.  In Section \ref{sec:Torelli} we prove Theorem \ref{thm:Torelli} and describe the relationship between Kirwan surjectivity and the action of the finite group $\Gamma_3$ on the cohomology of the space of semistable points.

\section{Stratifications}\label{sec:stratifications}

\subsection{Critical points of the Yang-Mills-Higgs functional}

The goal of this section is to describe the stratification of the space $\Bcal(d_1,d_2)$. There are in fact two natural stratifications: the Morse stratification given by the gradient flow of the Yang-Mills-Higgs functional which is detailed in this subsection, and the algebraic stratification according to Harder-Narasimhan type which is discussed in the next subsection. In Proposition \ref{prop:stratifications} we claim that, as in \cite{DWWW11} and \cite{WentworthWilkin11}, these stratifications coincide. 

We begin with the  classification of the critical sets of the Yang-Mills-Higgs functional in the
general case.  Fix
smooth complex hermitian vector bundles $E_p, E_q\to X$, with $\rank E_p=p$, $\rank E_q=q$, $\deg E_i=d_i$.
Without loss of generality (see the remark in \cite[p731]{Gothen02}) we always assume that $pd_q \geq qd_p$.
A \emph{$\U(p,q)$ Higgs bundle} consists of a split holomorphic structure on $V = 
E_p \oplus E_q$, and a Higgs field of the type $H^0(E_p^* E_q \otimes K) \oplus H^0(E_q^* E_p\otimes
K)$, where $K$ is the canonical bundle of $X$. In other words, a pair
\begin{equation*}
(\bar{\partial}_A, \Phi) = \left( \bar{\partial}_{A_p} \oplus \bar{\partial}_{A_q} , \left( \begin{matrix} 0 & c \\ b & 0 \end{matrix} \right) \right) 
\end{equation*}
Let  $\Bcal(E_p,E_q)$ denote the space of $\U(p,q)$ Higgs bundles (for more details, see
\cite{BGG03}).

Equivalently, one may view 
a $\U(p, q)$ Higgs bundle as a \emph{twisted quiver bundle}, with diagram
\begin{equation} \label{eqn:quiver}
\xygraph{
!{<0cm, 0cm>;<1.0cm, 0cm>:<0cm, 1.0cm>::}
!{(0,0) }*+{\bullet_{E_q}}="a"
!{(2,0) }*+{\bullet_{E_p}}="b"
"a" :@/^0.3cm/ "b"  ^{c}
"b" :@/^0.3cm/ "a"  ^{b}}
\end{equation}
where $b : E_p \rightarrow E_q \otimes K$ and $c : E_q \rightarrow E_p \otimes K$ are holomorphic.
The Yang-Mills-Higgs functional on the space of Higgs bundles is defined by
\begin{equation*}
\YMH(\bar{\partial}_A, \Phi) = \left\| F_A + [\Phi, \Phi^*] \right\|^2 
\end{equation*}
where $F_A$ denotes the curvature of the Chern connection associated to $\dbar_A$
and the hermitian structure on $E_p\oplus E_q$,
and $\Vert\cdot\Vert$ is the $L^2$-norm taken with  respect to 
 a choice of conformal metric on $X$.
On restriction to the space $\B(E_p,E_q)$ this becomes
\begin{equation*}
\YMH(\bar{\partial}_{A_p}, \bar{\partial}_{A_q}, b, c) = \left\| F_{A_p} + b b^* + c^* c \right\|^2 + \left\| F_{A_q} + b^* b + c c^* \right\|^2 
\end{equation*}

Let $\Bcal_{min}(E_p,E_q)$ denote the set of absolute minima of $\YMH$.  More generally, 
the critical point equations for $\YMH$ on  $\mathcal{B}(E_p,E_q)$ are
\begin{align}
\bar{\partial}_{A_q} * \left( F_{A_q} + b b^* + c^* c \right) & = 0 \label{eqn:bundle-1} \\
\bar{\partial}_{A_p} * \left( F_{A_p} + b^* b + c c^* \right) & = 0 \label{eqn:bundle-2} \\
b * \left( F_{A_p} + b^* b + c c^* \right) - * \left( F_{A_q} + b b^* + c^* c \right) b & = 0 \label{eqn:Higgs-b} \\
c * \left( F_{A_q} + b b ^* + c^* c \right) - * \left( F_{A_p} + b^* b + c c ^* \right) c & = 0 \label{eqn:Higgs-c}
\end{align}
Using \eqref{eqn:bundle-1} and \eqref{eqn:bundle-2} and 
the same method of proof  for holomorphic bundles in \cite[Section 5]{AtiyahBott83}, we conclude
 that the eigenvalues of $* \left( F_{A_q} + b b^* + c^* c \right)$ and $* \left( F_{A_p} + b^* b + c c^* \right)$ are  constant and  the holomorphic structures on $E_p$ and $E_q$ split according to these eigenvalues. We can therefore write $* \left( F_{A_q} + b b^* + c^* c \right)$ and $* \left( F_{A_p} + b^* b + c c^* \right)$ in the following block-diagonal form
\begin{equation*}
*(F_{A_q} + b b^* + c^* c) = \left( \begin{matrix} \lambda_1^q & 0 & \cdots & 0 \\ 0 & \lambda_2^q & \cdots & 0 \\ \vdots & \vdots & \ddots & \vdots \\ 0 & 0 & \cdots & \lambda_{n_q}^q \end{matrix} \right)  
\quad \text{and} \quad * \left( F_{A_p} + b^* b + c c^* \right) = \left( \begin{matrix} \lambda_1^p & 0 & 
\cdots & 0 \\ 0 & \lambda_2^p & \cdots & 0 \\ \vdots & \vdots & \ddots & \vdots \\ 0 & 0 & 
\cdots & \lambda_{n_p}^p \end{matrix} \right) 
\end{equation*}
(recall that these expressions are skew-Hermitian with respect to the metrics on $E_p$ and $E_q$, and hence diagonalizable). The bundles $E_p$ and $E_q$ then split with respect to this decomposition as follows
$$
E_p  = E_p^{(\lambda_1^p)} \oplus \cdots \oplus E_p^{(\lambda_{n_p}^p)} \ , \
E_q  = E_q^{(\lambda_1^q)} \oplus \cdots \oplus E_q^{(\lambda_{n_q}^q)} 
$$
where $E_p^{(\lambda_k^p)}$ (resp.\ $E_q^{(\lambda_k^q)}$) is the holomorphic sub-bundle of $E_p$
(resp.\ $E_q$) corresponding to the eigenvalue $\lambda_k^p$ (resp. $\lambda_k^q$). 
The Higgs fields $b$ and $c$ also decompose with respect to this splitting, 
and it follows from equations 
\eqref{eqn:Higgs-b} and \eqref{eqn:Higgs-c} that, if $\lambda_{j}^p \neq \lambda_{k}^q$, 
then the component of $b$ mapping $E_p^{(\lambda_{j}^p)}$ to $E_q^{(\lambda_{k}^q)}$ 
is zero and the component of $c$ mapping $E_q^{(\lambda_{k}^q)}$ to $E_p^{(\lambda_{j}^q)}$ is zero.

Therefore, the critical point equations define a splitting of $(\bar{\partial}_{A_p}, 
\bar{\partial}_{A_q}, b, c)$ into $\U(p', q')$ sub-bundles
\begin{equation}
(\bar{\partial}_{A_p}, \bar{\partial}_{A_q}, b, c) = \bigoplus_{\ell} (\bar{\partial}_{A_p^\ell}, \bar{\partial}_{A_q^\ell}, b_\ell, c_\ell) 
\end{equation}
where $\ell$ ranges over the set of all eigenvalues of 
$* \left( F_{A_q} + b b^* + c^* c \right)$ and 
$* \left( F_{A_p} + b^* b + c c^* \right)$, and $q' = \rank(E_q^{(\ell)})$, 
$p' = \rank(E_p^{(\ell)})$ 
(note that it is possible for one of $p'$ or $q'$ to be zero). Moreover, the usual Chern-Weil technique shows that the eigenvalues are determined by the slope of the bundles $E_p^\ell \oplus E_q^\ell$, and that $(\bar{\partial}_{A_q^\ell}, \bar{\partial}_{A_p^\ell}, b_\ell, c_\ell)$ minimizes the Yang-Mills-Higgs functional on $\mathcal{B}(E_p^\ell, E_q^\ell)$.

These results are summarized in the following proposition.

\begin{proposition}
A $\U(p,q)$ Higgs structure $(\bar{\partial}_{A_p}, \bar{\partial}_{A_q}, b, c)$ is a critical point for the Yang-Mills-Higgs functional if and only if it splits into the direct sum of $U(p', q')$ sub-bundles, each of which is a minimizer for the associated Yang-Mills-Higgs functional on the sub-bundles. The splitting is determined by the eigenvalues and eigenspaces of $* \left( F_{A_q} + b b^* + c^* c \right)$ and $* \left( F_{A_p} + b^* b + c c^* \right)$.
\end{proposition}

For the special case $p=2$ and $q=1$, there are only three types of decomposition that can occur at nonminimal critical points. 
The first is where the Higgs field is zero and the bundle $E_2$ is polystable.
 This corresponds to a splitting of the structure into a $\U(1)$ and a $\U(2)$ Higgs bundle.
 Call these critical points \emph{Type A} and let $\mathcal{C}_a$ denote the set of all critical points of Type A.

The second type of decomposition is where the Higgs field is zero and the structure splits into three $\U(1)$ Higgs  bundles. Call these critical points \emph{Type B}. In this case, the bundle $E_2$ is the direct sum of holomorphic line bundles $E_2 \cong S \oplus Q$. Without loss of generality, assume that $\deg S > \deg Q$, and note that the Higgs field is necessarily zero when $\deg S \neq d_1$.  
Also, use the notation $d_S=\deg S$, $d_Q=\deg Q$.

For convenience, when $d_S = d_1$ we also include the possibility that the Higgs field can be nonzero (see also Remark \ref{rem:exceptional-cases}). The critical point equations imply that such a Higgs bundle must take the following form
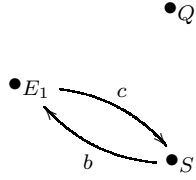
\begin{equation*}
\xygraph{
!{<0cm, 0cm>;<1.0cm, 0cm>:<0cm, 1.0cm>::}
!{(2,0) }*+{\bullet_{S}}="a"
!{(0,1) }*+{\bullet_{E_1}}="b"
!{(2,2) }*+{\bullet_{Q}}="c"
"b" :@/^0.3cm/ "a"  ^{c} "a" :@/^0.3cm/ "b" ^{b} }
\end{equation*}
where $b$ and $c$ are related by $\| b \|^2 = \| c \|^2$. 

The connected components of the space of Type B critical points are in 
one-to-one correspondence with the range of values for $d_S$.
 Moreover, there are three different cases for $d_S$ 
that lead to different contributions to the Morse theory calculations of
 Section \ref{sec:calculations}. 
For each value of $\ell$ in the range $\frac{1}{2} d_2 < \ell < d_1$,
let $\mathcal{C}_{b_1}^\ell$ denote the set of Type B critical
points for which $d_S = \ell$, 
define $\mathcal{C}_{b_2}^{d_1}$ to be the set of Type B critical
points for which $d_S = d_1$, 
and for $d_1 < \ell$ define $\mathcal{C}_{b_3}^\ell$ to be the set of
Type B critical points for which $d_S = \ell$.

The third type of decomposition is where the $\U(2,1)$ structure splits 
into the direct sum of a stable $U(1,1)$ structure and a $U(1)$ structure. 
Equivalently, the bundle $E_2$ splits into line bundles, 
$E_2 \cong S \oplus Q$, and, depending on the degree of $S$ and $Q$, 
the Higgs field takes on one of the following forms.
\begin{enumerate}
\item $d_S > \frac{1}{2}(d_Q + d_1)$.
 In this case the maximal semistable subobject of the Higgs bundle
$(E_2 \oplus E_1, b, c)$ is a line subbundle of $S$, which does not
interact with the Higgs field. Define  $\ell = d_S$. Since we have assumed that $d_2 \leq
2d_1$ (see the Introduction), then the condition $\ell > \frac{1}{2} (d_Q + d_1)$ implies
that $d_1 > d_2 - \ell = d_Q$. Minimality of the Yang-Mills-Higgs
functional on the subobject $(Q \oplus E_1, b,c)$ then implies that $b$
and $c$ are related by $\frac{1}{\pi}\left( \| c \|^2 - \| b
\|^2 \right) = d_1 - d_Q > 0$
 and therefore $c \neq 0$. Label these critical sets $\mathcal{C}_{c_1}^{\ell}$. A graphical representation of the Higgs field at these critical points is as follows.
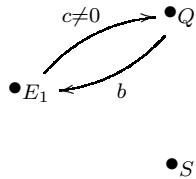
\begin{equation*}
\xygraph{
!{<0cm, 0cm>;<1.0cm, 0cm>:<0cm, 1.0cm>::}
!{(2,0) }*+{\bullet_{S}}="a"
!{(0,1) }*+{\bullet_{E_1}}="b"
!{(2,2) }*+{\bullet_{Q}}="c"
"b" :@/^0.3cm/ "c"  ^{c \neq 0} "c" :@/^0.3cm/ "b" ^{b} }
\end{equation*}
The section $c$ can only be nonzero if $\deg (E_1^* Q \otimes K) \geq 0$, and so these critical points only exist for values of $\ell$ such that
$
d_2 - \ell - d_1 + 2g-2 \geq 0
$ and $\ell > \tfrac{1}{2} (d_2 - \ell + d_1)$.
This is equivalent to the condition that $\ell$ is in the range
$
\frac{1}{3}(d_1 + d_2) < \ell \leq d_2 - d_1 + 2g-2  
$.

\item $d_Q < \frac{1}{2}(d_S + d_1)$ and $d_1 > d_S$. In this case the maximal semistable subobject of $(E_1, E_2, b, c)$ is 
$(E_1 \oplus S, b, c)$, and so we define $\ell = d_S$ 
and $ d_Q=d_2-\ell$. Then the same analysis as before shows that $b$
and $c$ are related by $\frac{1}{\pi}\left( \| c \|^2 - \| b
\|^2 \right) = d_1 - d_S > 0$, and therefore $c \neq 0$. Call these critical sets $\mathcal{C}_{c_2}^{\ell}$. The corresponding picture is
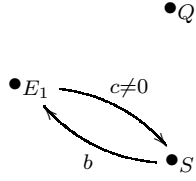
\begin{equation*}
\xygraph{
!{<0cm, 0cm>;<1.0cm, 0cm>:<0cm, 1.0cm>::}
!{(2,0) }*+{\bullet_{S}}="a"
!{(0,1) }*+{\bullet_{E_1}}="b"
!{(2,2) }*+{\bullet_{Q}}="c"
"a" :@/^0.3cm/ "b"  ^{b} "b" :@/^0.3cm/ "a" ^{c \neq 0} }
\end{equation*}
Critical sets of this type can only exist if $c \neq 0$, and so we must have $\deg (E_1^*S\otimes K) \geq 0$. Combining this with the conditions that $d_Q < \frac{1}{2} (d_S + d_1)$ and $d_1 > d_S$ gives 
\begin{equation}\label{eqn:prelim-c2-range}
\max \bigl( d_1 - 2g+1, \tfrac{1}{3} (2d_2 - d_1) \bigr) < \ell < d_1 
\end{equation}
The bound on the Toledo invariant $2d_1 -d_2 \leq 3g-3$  is equivalent to $d_1 - (2g-2) \leq \frac{1}{3}(2d_2-d_1)$. Therefore, the inequality \eqref{eqn:prelim-c2-range} reduces to
$
\tfrac{1}{3} (2d_2-d_1) < \ell < d_1$.

\item $d_Q < \frac{1}{2}(d_S + d_1)$ and $d_1 < d_S$. In this case
the maximal semistable subobject of $(E_2 \oplus E_1, b, c)$ is
$(E_1 \oplus S, b, c)$, and so we define $\ell = d_S$ and $d_Q=d_2-\ell$. 
An analysis of the critical point equations shows that now $b \neq 0$
and that $b$ and $c$ are related by $\frac{1}{\pi}\left( \| b
\|_{L^2}^2 - \| c \|_{L^2}^2 \right) = d_S - d_1 > 0$. Call these critical sets $\mathcal{C}_{c_3}^{\ell}$, and note that the quiver bundle picture reduces to
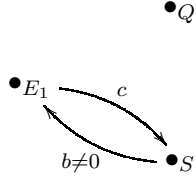
\begin{equation*}
\xygraph{
!{<0cm, 0cm>;<1.0cm, 0cm>:<0cm, 1.0cm>::}
!{(2,0) }*+{\bullet_{S}}="a"
!{(0,1) }*+{\bullet_{E_1}}="b"
!{(2,2) }*+{\bullet_{Q}}="c"
"a" :@/^0.3cm/ "b"  ^{b \neq 0} "b" :@/^0.3cm/ "a" ^{c} }
\end{equation*}
These critical sets can only exist if $b \neq 0$, and so $\deg (S^* E_1 \otimes K) \geq 0$. Note that $d_1 < \ell$ implies that $d_Q < \frac{1}{2}(d_S + d_1)$, and so we have the inequalities
$
d_1 - \ell + 2g-2 \geq 0$
and $\ell > d_1$.
Therefore
$d_1 < \ell \leq d_1 + 2g-2$.

\end{enumerate}

\begin{remark}\label{rem:exceptional-cases}
Note that there are two possible values of $\ell$ that have not been included in the above list. The first is $\ell = d_1$, for which the critical points have already been classified as type $\mathcal{C}_{b_2}^{d_1}$. The second is $\ell = \frac{1}{3}(2d_2 - d_1)$, in which case the critical point minimizes the Yang-Mills-Higgs functional.
\end{remark}


Using the descriptions above, a standard calculation gives the following table of results for the equivariant Poincar\'e polynomial of each nonminimal critical set.
\begin{table}[h]
\caption{Classification of the critical sets and their topology}
\centering

\begin{tabular}{ c | c | c }
\hline \hline
Critical set & Range of values of $\ell$ & Equivariant Poincar\'e polynomial \\
[1ex] \hline 
 $\mathcal{C}_a$ & n/a & $\frac{1}{(1-t^2)^2} P_t(J(X)) P_t^{\Gcal(E_2)} (\mathcal{A}^{ss}(E_2))$  \\ [1ex] 
$\mathcal{C}_{b_1}^{\ell}$ & $\frac{1}{2} d_2 < \ell < d_1$ & $\frac{1}{(1-t^2)^3} P_t(J(X))^3$  \\ [1ex] 
$\mathcal{C}_{b_2}^{\ell}$ & $\ell = d_1$ & $\frac{1}{(1-t^2)^3} P_t(J(X))^3$  \\ [1ex] 
$\mathcal{C}_{b_3}^{\ell}$ & $d_1 < \ell$ & $\frac{1}{(1-t^2)^3} P_t(J(X))^3$ \\ [1ex] 
$\mathcal{C}_{c_1}^{\ell}$ & $ \frac{1}{3}(d_1 +d_2)< \ell \leq d_2-d_1 +2g-2$ & $\frac{1}{(1-t^2)^2} P_t(J(X))^2 P_t(S^{\ell - d_1 + 2g-2} X)$ \\ [1ex] 
$\mathcal{C}_{c_2}^{\ell}$ & $\frac{1}{3}(2d_2 - d_1) < \ell < d_1$ & $\frac{1}{(1-t^2)^2} P_t(J(X))^2 P_t(S^{\ell - d_1 + 2g-2} X)$ \\ [1ex] 
$\mathcal{C}_{c_3}^{\ell}$ & $d_1 < \ell \leq d_1+ 2g-2$ & $\frac{1}{(1-t^2)^2} P_t(J(X))^2P_t(S^{d_1 - \ell + 2g-2} X)$ \\ [1ex] 
\hline
\end{tabular}
\end{table}
In  Table 1 we have used the following notation
 for the Harder-Narasimhan stratification for holomorphic bundles. 
Suppose $E\to X$ is a  rank 2 bundle of degree $d$. 
Let $\Acal(E)$ denote the space of $\dbar$-operators on $E$, which 
is equivalent to the space of holomorphic structures. 
 Let $\Acal^{ss}(E)\subset \Acal(E)$ denote the subset of semistable bundles, and 
  for $j>d/2$, let $\Acal_j(E)\subset\Acal(E)$ be the subset of 
unstable holomorphic structures on  $E$ whose maximally destabilizing line subbundle has degree $j$.

We denote the ordered set of possible values in the labeling of the critical sets above by
\begin{equation}  \label{eqn:delta}
\Delta_{d_1,d_2}=\{ \tfrac{1}{2}d_2\}\cup
 \left\{ \ell\in \ZBbb : \ell>
\tfrac{1}{3}(2d_2-d_1) \right\} 
\end{equation}
 We will express the
 various components as $\Ccal_a$, $\Ccal_{b}^\ell$, and $\Ccal_c^\ell$.

\subsection{Harder-Narasimhan and Morse stratifications}\label{sec:H-N-and-Morse}
We now describe the algebraic stratification of the space of $\U(2,1)$ Higgs bundles. As in the previous section, let
$$
V=E_2\oplus E_1\ ,\ \Phi=\left(\begin{matrix} 0& c\\ b&0\end{matrix}\right)
$$
Recall that $(V,\Phi)$ is \emph{stable} (resp.\ \emph{semistable}) if
$$
\mu(F)=\frac{\deg F}{\rank F}<\mu(V)=\frac{\deg V}{\rank V}\qquad\text{(resp.\ $\leq$)}
$$
for every $\Phi$-invariant subsheaf $0\neq\rank F\neq\rank V$.
If $(V,\Phi)$ is not semistable, a \emph{maximally destabilizing subbundle} is a subsheaf $0\neq F\subsetneq V$,  satisfying the following:
\begin{itemize}
\item $F$ is $\Phi$-invariant;
\item $\mu(F)> \mu(V)$;
\item $F$ is maximal in the sense that for any $F'\neq F$ satisfying the first two conditions, then $\mu(F')\leq \mu(F)$, and if equality, then $\rank F'<\rank F$.
\end{itemize}
If $F$ satisfies these conditions then $F$ must be saturated, i.e.\ $V/F$ is torsion-free.  
Unstable Higgs bundles have a unique (Harder-Narasimhan) filtration by sub-Higgs bundles.  The associated graded of this filtration will be denoted by $\Gr(V,\Phi)$.

Recall that by assumption,  $2d_1\geq d_2$. Below we  determine all the possible Harder-Narasimhan filtrations of unstable $\U(2,1)$ Higgs bundles. Let $F$ be a maximally destabilizing subbundle of $(V,\Phi)$.

\medskip
\noindent  {\bf Case I:} $\rank F=1$.
Let $f_i$ be the induced maps $F\to E_i$. 
 Then $f_2\equiv 0$ implies $f_1$ is an isomorphism, and $f_1\equiv 0$ implies $f_2$ is everywhere injective, and we claim that one of these two possibilities occurs.  For
suppose  neither $f_i\equiv 0$,  and let $F_2\subset E_2$ be the saturation of $\im
f_2$.
Then $E_1 \oplus F_2$
is a  subbundle  with slope at least $\deg F$, contradicting the assumption that $F$ is maximal.
It follows that there are two possibilities according to whether $F$ lies in $E_1$ or $E_2$.
\begin{enumerate}
\item If $F = E_1$, then $c \equiv 0$. If $E_2$ is semistable then the stratum is defined by the condition $c \equiv 0$, and we label it by $\mathcal{S}_a$. The quiver diagram in this case is
\begin{equation*}
\xygraph{
!{<0cm, 0cm>;<1.0cm, 0cm>:<0cm, 1.0cm>::}
!{(0,0) }*+{\bullet_{E_1}}="a"
!{(2,0) }*+{\bullet_{E_2}}="b"
"b" !{\ar @{-->}_{\ b} "a"}}
\end{equation*}
and the associated graded is $(E_2,0)\oplus (E_1,0)$
(In this diagram and the others below, we use a \emph{dashed arrow} to represent
 a component of the Higgs field that may or may not be zero and a 
\emph{solid arrow} to represent a component of the Higgs field that must be nonzero. 
If a component of the Higgs field must be zero then there is no arrow 
between the vertices). If the bundle $E_2$ is unstable, let $S\subset E_2$ be the
maximal destabilizing line bundle, and write  
 $0 \rightarrow S \rightarrow E_2 \rightarrow Q \rightarrow 0$, 
 with extension class $[a]\in H^1(X,Q^\ast S)$.
Notice that $d_S=\deg S < d_1$, 
since either $S \subset \ker b$ and $S$ 
is a subobject of $(V,\Phi)$, 
or $S$ is not in $\ker b$ and $S 
\oplus E_1$ is a subobject. 
The associated graded $\Gr(V,\Phi)=(E_1, 0)\oplus (S, 0)\oplus (Q,0)$, and we 
label this stratum by $\mathcal{S}_{b_1}^{\ell}$, where $\ell = \deg S$.
The quiver diagram for this case is
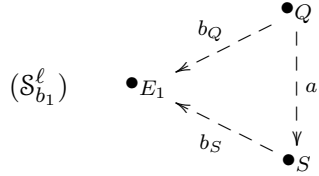
\begin{equation*}
 (\mathcal{S}_{b_1}^{\ell})\   \hskip.5cm
 \xygraph{
!{<0cm, 0cm>;<1.0cm, 0cm>:<0cm, 1.0cm>::}
!{(2,-1) }*+{\bullet_{S}}="a"
!{(0,0) }*+{\bullet_{E_1}}="b"
!{(2,1) }*+{\bullet_{Q}}="c"
"a" !{\ar @{-->}^{b_S} "b"} "c" !{\ar @{-->}_{b_Q} "b"}
"c" !{\ar @{-->}^{a} "a"}  }
\end{equation*}

\item If $F= S \subset E_2$ with quotient $Q$, then $S \subset \ker b$, 
$E_2$ is unstable, and the graded object of the Harder-Narasimhan filtration of 
$E_2$ is precisely $S \oplus Q$.  If 
$c_Q\equiv 0$, then we also require $d_S > d_1$, for otherwise $E_1$ would
 be invariant with slope at least $d_S$. In this case,
  $\Gr(V, \Phi)=(S, 0)\oplus (E_1,0)\oplus (Q,0)$.
   If 
$c_Q\neq 0$, then the only requirement is that $d_S > \frac{1}{3}(d_1 + d_2)$ 
(otherwise $E_1 \oplus Q$ would be invariant with slope at least $d_S$), and
$\Gr(V,\Phi)=(S, 0)\oplus   ( E_1 \oplus Q, 
b_Q, c_Q)$, where $(b_Q,c_Q)$ is the induced Higgs field on $E_1 \oplus Q$ coming 
from $b$ and the projection $c_Q$ of $c$ to $Q$. 
We label the strata $\mathcal{S}_{b_3}^{\ell}$ and
$\mathcal{S}_{c_1}^{\ell}$, respectively. The quiver diagrams for the two cases are
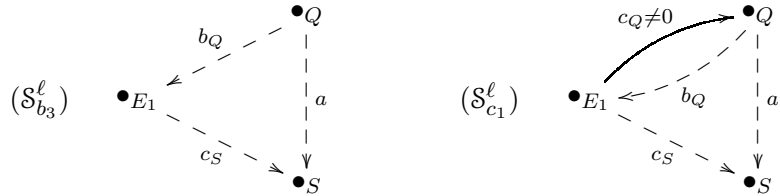
\begin{equation*}
(\mathcal{S}_{b_3}^{\ell}) \hskip.5cm \xygraph{
!{<0cm, 0cm>;<1.5cm, 0cm>:<0cm, 1.5cm>::}
!{(1.5,-.75) }*+{\bullet_{S}}="a"
!{(0,0) }*+{\bullet_{E_1}}="b"
!{(1.5,.75) }*+{\bullet_{Q}}="c"
"c" !{\ar @{-->}_{b_Q} "b"} "b" !{\ar @{-->}_{c_S} "a"}  
"c" !{\ar @{-->}^{a} "a"}  
!{(3.25,0) }*+{ (\mathcal{S}_{c_1}^{\ell})\quad }  
!{(5.5,-.75) }*+{\bullet_{S}}="a"
!{(4,0) }*+{\bullet_{E_1}}="b"
!{(5.5,.75) }*+{\bullet_{Q}}="c"
"b" :@/^0.3cm/ "c"  ^{c_Q \neq 0}
 "c" :@/^.3cm/  @{-->} "b" ^{b_Q}
"b" !{\ar @{-->}_{c_S} "a"}
"c" !{\ar @{-->}^{a} "a"}  }  
\end{equation*}
\end{enumerate}

\medskip
\noindent  {\bf Case II:} $\rank F=2$.  The projection $F\to E_1$ cannot vanish.
  Indeed, if if did, then $F=E_2$ and $d_2/2> (1/3)(d_1+d_2)$. But this 
contradicts the assumption $d_2\leq 2d_1$.   Let $S$ be the kernel of the 
projection $F\to E_1$.  Then $\deg(P=F/S) \leq d_1$.
We also have $S \subset E_2$.  Since $E_2/S$ is a subsheaf of $V/F$
 which we assume to be torsion-free, we conclude that $S$ is
 a subbundle of $E_2$. 
Let $[a_F]$ and $[a]$ denote the extension classes for the sequences
\begin{align}
0\lra S \lra F \lra P \lra 0 \label{eqn:extF} \\
0\lra S \lra E_2 \lra Q \lra 0 \label{eqn:ext2} 
\end{align}
In terms of the smooth splittings $E_2 \oplus E_1 = S \oplus Q\oplus E_1$ and $F = S \oplus P$, we can write the inclusion $F\hookrightarrow V$ and the Higgs field as
$$
f=\left(\begin{matrix}  1&f_1  \\ 0& f_2  \\ 0&f_P   \end{matrix}\right)\ ,\ 
\Phi=\left(\begin{matrix} 0 & 0 & c_S \\ 0 & 0 & c_Q \\ b_S & b_Q & 0\end{matrix}\right)
$$
where $f_P : P \rightarrow E_1$ is nonzero (since the projection of $F$ to $E_1$ cannot vanish), and $f_1 : P \rightarrow S$, $f_2 : P \rightarrow Q$ are induced by the projection from $F$ to $E_2$. 
Since $f$ has everywhere rank $2$,  $f_2$ and $f_P$ have no common zeros.
Holomorphicity of $f$ implies $f_2, f_P$ holomorphic, and $f_1$ satisfies
\begin{equation} \label{eqn:f1}
\dbar f_1+a f_2-a_F=0 
\end{equation}
where $\dbar$ is the induced holomorphic structure on $P^\ast \otimes S$. On the other hand, since $F$ is destabilizing
\begin{align*}
\deg(Q P^\ast) &= d_Q -\deg P = d_2-d_S -\deg P
=d_2-\deg F \\
&< d_2-\tfrac{2}{3}(d_1+d_2) 
= -\tfrac{1}{3}(2d_1-d_2)
\leq 0
\end{align*}
by the assumption on degrees.  It follows that $f_2\equiv 0$,  $f_P$ gives an isomorphism $P\cong E_1$, and by 
\eqref{eqn:f1} the sequence \eqref{eqn:extF} splits.
The condition that $F \cong E_1 \oplus S$ be
 invariant under the Higgs field is equivalent to $c_Q \equiv 0$. Moreover, $S$ is invariant if and only if $b_S \equiv 0$.
These are  the only conditions coming from invariance.

A splitting of $F\subset V$ gives a splitting of \eqref{eqn:ext2}. In this case, 
$\Gr(V,\Phi)=(E_1\oplus S, b_S, c_S)\oplus (Q,0)$.
and the condition on degrees is $\ell> \tfrac{1}{3}(2d_2-d_1)$. 
By the assumption that $F$ is maximally destabilizing, 
$\ell< d_1\Rightarrow  c_S \neq 0$, and 
$\ell > d_1  \Rightarrow  b_S \neq 0$.
We label the former case $\mathcal{S}_{c_2}^{\ell}$ and the latter case $\mathcal{S}_{c_3}^{\ell}$. 
When $\ell = d_1$ then there are no conditions on $b_S$ and $c_S$. We label this stratum $\mathcal{S}_{b_2}$.
The quiver diagrams for these strata are

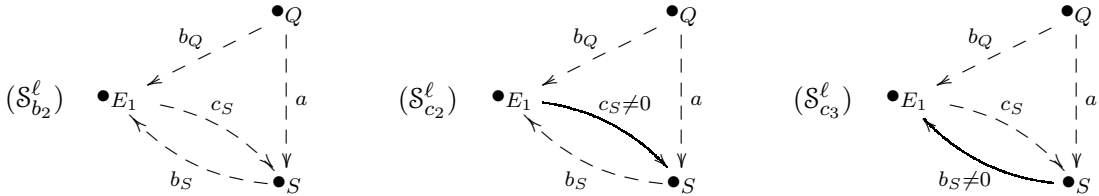
\begin{equation*}(\mathcal{S}_{b_2}^{\ell}) \hskip.3cm \xygraph{
!{<0cm, 0cm>;<1.5cm, 0cm>:<0cm, 1.5cm>::}
!{(1.5,-.75) }*+{\bullet_{S}}="a"
!{(0,0) }*+{\bullet_{E_1}}="b"
!{(1.5,.75) }*+{\bullet_{Q}}="c"
"c" !{\ar @{-->}^{a} "a"}  
"b" !{\ar @/^0.3cm/  @{-->}^{c_S} "a"} "c" !{\ar @{-->}_{b_Q} "b"} "a" !{\ar @/^0.3cm/ @{-->}^{b_S} "b"}  
!{(2.75,0) }*+{ (\mathcal{S}_{c_2}^{\ell})}  
!{(5,-.75) }*+{\bullet_{S}}="a"
!{(3.5,0) }*+{\bullet_{E_1}}="b"
!{(5,.75) }*+{\bullet_{Q}}="c"
"c" !{\ar @{-->}^{a} "a"}  
"b" :@/^0.3cm/ "a"  ^{c_S \neq 0}
"c" !{\ar @{-->}_{b_Q} "b"} "a" !{\ar @/^0.3cm/ @{-->}^{b_S} "b"}  
!{(6.25,0) }*+{ (\mathcal{S}_{c_3}^{\ell})}  
!{(8.5,-.75) }*+{\bullet_{S}}="a"
!{(7,0) }*+{\bullet_{E_1}}="b"
!{(8.5,.75) }*+{\bullet_{Q}}="c"
"a" :@/^0.3cm/ "b"  ^{b_S \neq 0} "c" !{\ar @{-->}_{b_Q} "b"} "b" !{\ar @/^0.3cm/ @{-->}^{c_S} "a"}  
"c" !{\ar @{-->}^{a} "a"} 
 }
 \end{equation*}

We conclude that there is a 1-1 correspondence between the associated graded objects 
listed above and  the critical sets of the
YMH functional.
The collection $\{\Scal_a, \Scal_b^\ell, \Scal_c^\ell\}$,  where
$\ell\in \Delta_{d_1,d_2}$ has its 
natural ordering, combine to form the  \emph{Harder-Narasimhan stratification} of $\Bcal(d_1,d_2)$.
As in \cite{WentworthWilkin11}, however, it turns out that this stratification is too fine a structure for an equivariantly perfect Morse theory. 
For $k\in \Delta_{d_1,d_2}$ define
\begin{align}
\label{eqn:HN-modified} X_k&=\begin{cases}\displaystyle  \Bcal^{ss}(d_1,d_2)\cup
\bigcup_{\ell\in \Delta_{d_1,d_2} ,\, \ell\leq k} \Scal_c^\ell   &\quad  k < d_2/2 \\
\displaystyle \Bcal^{ss}(d_1,d_2)\cup \Scal_a\cup\bigcup_{\ell\in \Delta_{d_1,d_2},\, \ell \leq d_2/2} \Scal_c^\ell   &\quad k=d_2/2 \\
\displaystyle X^{ss}\cup \Scal_a\cup\bigcup_{\ell\in \Delta_{d_1,d_2} ,\, \ell\leq k} \Scal_c^\ell  \cup \bigcup_{\ell\in \Delta_{d_1,d_2} ,\ \ell\leq k} \Scal_b^\ell  &\quad k > d_2/2
\end{cases} \\
\label{eqn:HN-star} X_k^\ast&=\begin{cases} \displaystyle \Bcal^{ss}(d_1,d_2)
\cup \bigcup_{\ell\in \Delta_{d_1,d_2} ,\, \ell < k} \Scal_c^\ell  &\quad k \leq d_2/2    \\
\displaystyle \Bcal^{ss}(d_1,d_2)
\cup \Scal_a\cup\bigcup_{\ell\in \Delta_{d_1,d_2} ,\, \ell < k} \Scal_c^\ell \cup \bigcup_{\ell\in \Delta_{d_1,d_2} ,\ \ell < k} \Scal_b^\ell  &\quad k > d_2/2
\end{cases}
\end{align}
 In the notation above, $\Scal^\ell_c$ means the union over all possible $C$-strata with index $\ell$.
This gives a  $\G$-invariant stratification of $\Bcal(d_1,d_2)$ which we refer to as the \emph{modified Harder-Narasimhan stratification}.

\begin{remark}
As in \cite{WentworthWilkin11}, the ordering of the set
$\Delta_{d_1,d_2}$ does not in general coincide with the one coming
from values of the YMH functional.  This is irrelevant for the
calculations in this paper.
\end{remark}

 As in
\cite{Wilkin08}, we have
\begin{proposition}  \label{prop:stratifications}
The Morse stratification of the YMH flow coincides with the
Harder-Narasimhan stratification of $\Bcal(d_1,d_2)$.  In particular, 
the gradient flow of the YMH functional defines $\Gcal$-equivariant deformation
retractions
 $\Bcal_{min}(d_1,d_2)\hookrightarrow \Bcal^{ss}(d_1,d_2)$,  $\Ccal_a\hookrightarrow \Scal_a$,
$\Ccal_b^\ell\hookrightarrow \Scal_b^\ell$, and 
$\Ccal_c^\ell\hookrightarrow \Scal_c^\ell$.
\end{proposition}

We now state one of the main results of this paper.  The proof occupies the next section.

\begin{theorem}[Perfect stratification] \label{thm:perfection}
The modified Harder-Narasimhan stratification $\{X_k\}_{k\in
\Delta_{d_1,d_2}}$ of $\Bcal(d_1,d_2)$  is $\Gcal$-equivariantly perfect
in the sense that the inclusions $\Bcal^{ss}(d_1,d_2)\subset X_{k}\subset X_\ell$ induce
surjections $H^\ast_\Gcal(X_k)\to H^\ast_\Gcal(\Bcal^{ss}(d_1,d_2))$ and $H^\ast_\Gcal(X_\ell)\to H^\ast_\Gcal(X_{k})$ for all
$k\leq \ell$ in $\Delta_{d_1,d_2}$.
\end{theorem}

\begin{corollary}[Kirwan surjectivity] \label{cor:kirwan}
The Kirwan  map $\kappa$  in \eqref{eqn:kirwan_map} is surjective.
\end{corollary}

In order to prove Theorem \ref{thm:perfection} we shall need to bootstrap an intermediary stratification lying between the HN and modified HN strata.  Define $\{X_k'\}_{k\in \Delta_{d_1,d_2}}$ by setting $X_k'=X_k^\ast\cup \Scal_c^k$.
There are three crucial regions,  essentially depending upon the number of $C$-strata.  We refer to these by the following:
\begin{enumerate}\label{def:regions}
\item[(${\bf I}$)]  where $\tfrac{1}{3}(d_1+d_2)<k\leq d_2-d_1+2g-2$;
\item[(${\bf II}$)] where $\tfrac{1}{3}(2d_2-d_1)<k\leq \tfrac{1}{3}(d_1+d_2)$, or where
$d_2-d_1+2g-2<k\leq d_1$ (if possible);
\item[(${\bf III}$)]  where $\max\{d_1,d_2-d_1+2g-2\}< k$. 
\end{enumerate}

\section{Singular Morse Theory}\label{sec:calculations}

In this section we develop the necessary machinery to perform the Morse theory calculations on the singular space $\mathcal{B}(E_1, E_2)$.
Recall Kirwan's result \cite{Kirwan84}. For a Hamiltonian action of a compact connected Lie group $K$ on a compact smooth symplectic manifold $M$, there is a compatible Riemannian structure such that
induced Morse stratification $\{ S_\mu \}_{\mu\in I}$ is smooth, where $I$ is a partially ordered set labeling the critical sets. Let 
$$X_\mu=\cup_{\nu\leq \mu} S_\nu\ ,\ X_\mu^\ast=\cup_{\nu< \mu} S_\nu
$$
Then Kirwan shows that the long exact sequence 
\begin{equation}\label{eqn:Kirwan-LES}
\cdots\lra H_K^p(X_\mu, X_\mu^\ast) \stackrel{\alpha^p}{\longrightarrow} H_K^p(X_\mu) \stackrel{\beta^p}{\longrightarrow} H_K^p(X_\mu^\ast) \longrightarrow \cdots 
\end{equation}
splits into short exact sequences. Moreover, the Thom isomorphism  implies that $H_K^p(X_\mu, X_\mu^\ast) \cong H_K^{p-\lambda_\mu}(C_\mu)$, where  $C_\mu$ is the critical set at the minimum of the stratum $S_\mu$ and $\lambda_\mu$ is the Morse index. The splitting of \eqref{eqn:Kirwan-LES} is a consequence of the fact that $\alpha^p$ is always injective, which in turn follows from the Atiyah-Bott lemma \cite{AtiyahBott83}. Therefore, to compute the change in cohomology that occurs when attaching the stratum $S_\mu$, it is sufficient to know the cohomology and the Morse index of each critical set. Moreover, $\alpha^p$ injective for all $p$ implies that $\beta^p$ is surjective for all $p$, and so inclusion $X_{\mu}^\ast \hookrightarrow X_\mu$ induces a surjective map $H_K^*(X_\mu) \rightarrow H_K^*(X_\mu^\ast)$. 

When the ambient space is singular, the idea behind the calculation is an extension of the one described above. We still study the long exact sequence \eqref{eqn:Kirwan-LES}, however the calculation of $H_K^p(X_\mu, X_\mu^\ast)$ is much more complicated than an application of the Thom isomorphism, and in fact $\alpha^p$ is not always injective for $\SU(2,1)$ Higgs bundles.

This section is divided into four subsections. 
In the first subsection we describe the negative eigenspace of the Hessian at each 
critical point. In the second we compute the relevant cohomology groups needed to compute 
$H_\mathcal{G}^*(\nu_\mu^-, \nu^-_\mu\setminus\{0\})$ (where $\nu^-_\mu$ denotes 
the negative normal space to the critical set). 
In the third subsection we prove the isomorphism $H_\mathcal{G}^*(X_\mu, X_\mu^\ast)
 \cong H_\mathcal{G}^*(\nu^-_\mu, \nu_\mu^-\setminus\{0\})$ (in certain cases), 
and in the final section we show that the modified Harder-Narasimhan stratification defined in the previous section is equivariantly perfect for $\U(2,1)$ Higgs bundles (i.e.\ our analog of \eqref{eqn:Kirwan-LES} splits into short exact sequences). 

Finally, it is worth mentioning here that \emph{a priori} we should do our cohomology computations on a 
small neighborhood of the zero section in the negative normal space $\nu^-_\mu$. 
The proofs of the relevant results 
(e.g.\ Proposition \ref{prop:bott1})  decompose all of the necessary calculations to calculations where the ambient space is a manifold, which allows us to study the whole space $\nu_\mu^-$ instead of a neighborhood of the zero section. This observation simplifies some of the definitions and calculations in this section.

\subsection{Indices of critical sets} \label{sec:index}

First, recall the following result for Higgs vector bundles.

\begin{lemma}
Let $(A, \Phi)$ be a critical point of $\YMH$ on the space of Higgs bundles. 
A pair $( \alpha,  \varphi) \in \Omega^{0,1}(\End V) \oplus \Omega^{1,0}(\End V)$ is in the negative eigenspace of the Hessian at $(A, \Phi)$ if
\begin{enumerate}

\item The pair $(\alpha,  \varphi)$ is orthogonal
 to the $\mathcal{G}^\CBbb$-orbit 
through $(A, \Phi)$. Equivalently, 
$\bar{\partial}_A^* \alpha - \bar{*} [ \Phi, \bar{*}  \varphi ] = 0
$.

\item The pair $(A + \alpha, \Phi +  \varphi)$ is a Higgs pair. Equivalently, the following equation is satisfied
\begin{equation*}
\bar{\partial}_A  \varphi + [ \alpha, \Phi] + [\alpha,  \varphi] = 0
\end{equation*}

\item The pair $(\alpha,  \varphi)$ is an eigenvector for the operator $i \ad *(F_A + [\Phi, \Phi^*])$ with negative eigenvalue. Equivalently, the following equations are satisfied
\begin{align*}
i\left[ *( F_A + [\Phi, \Phi^*] ), \alpha \right] & = \lambda \alpha \\
i\left[ *( F_A + [\Phi, \Phi^*] ),  \varphi \right] & = \lambda  \varphi
\end{align*}
for some $\lambda < 0$. (Note that the eigenvalues are necessarily real since $i *(F_A + [\Phi, \Phi^*])$ is self-adjoint.)

\end{enumerate}

\end{lemma}

To translate this into a statement for $\U(p, q)$ Higgs bundles $V=E_p\oplus E_q$,  we use the following inclusions
\begin{align*}
\Omega^{0,1}(\End E_p) \oplus \Omega^{0,1}(\End E_q) &\hookrightarrow \Omega^{0,1}(\End V) \\ 
\Omega^0(E_p^* E_q \otimes K) \oplus \Omega^0(E_q^* E_p \otimes K) &\hookrightarrow \Omega^0((\End V) \otimes K) 
\end{align*}

\begin{corollary}
Let $(A_p, A_q, b, c)$ be a critical point of $\YMH$ on
$\Bcal(E_p,E_q)$. Then 
\begin{equation*}
(\alpha_p, \alpha_q, \beta, \gamma) \in \Omega^{0,1}(\End E_p ) \oplus \Omega^{0,1}(\End E_q) \oplus \Omega^0(E_p^* E_q \otimes K) \oplus \Omega^0(E_q^* E_p \otimes K)
\end{equation*}
is in the negative eigenspace of the Hessian at $(A_p, A_q, b, c)$ if
\begin{enumerate}

\item $(\alpha_p, \alpha_q, \beta, \gamma)$ is orthogonal to the
 $\mathcal{G}^\CBbb$-orbit through $(A_p, A_q, b, c)$. 
Equivalently, the following equations are satisfied
\begin{align*}
\bar{\partial}_{A_q}^* \alpha_q  - \bar{*} (b (\bar{*} \beta)) - \bar{*} ((\bar{*} \gamma) c) & = 0 \\
\bar{\partial}_{A_p}^* \alpha_p - \bar{*} ( (\bar{*} \beta) b) - \bar{*} ( c (\bar{*} \gamma)) & = 0 
\end{align*}

\item $(A_p + \alpha_p, A_q + \alpha_q, b + \beta, c + \gamma)$ is a Higgs pair. Equivalently, the following equations are satisfied
\begin{align*}
\bar{\partial}_A \beta + (\alpha_q) (b + \beta) + (b + \beta) (\alpha_p) & = 0 \\
\bar{\partial}_A \gamma + (\alpha_p) (c + \gamma) + (c + \gamma) (\alpha_q) & = 0 
\end{align*}
where $\bar{\partial}_A$ denotes the holomorphic structure induced by $\bar{\partial}_{A_p}$ and 
$\bar{\partial}_{A_q}$ on both $E_p^* E_q \otimes K$ and $E_q^* E_p \otimes K$.

\item The pair $(\alpha,  \varphi)$ is an eigenvector for the operator $i \ad *(F_A + [\Phi, \Phi^*])$ with negative eigenvalue. Equivalently, the following equations are satisfied
\begin{align*}
i\left[ *( F_{A_q} + b b^* + c^* c ), \alpha_q \right] & = \lambda \alpha_q \\
i\left[ *( F_{A_p} + b^* b + c c^* ), \alpha_p \right] & = \lambda \alpha_p \\
i \left( *( F_{A_q} + b b^* + c^* c) \beta - \beta  *( F_{A_p} + b^* b + c c^* ) \right) & = \lambda \beta \\
i \left( *( F_{A_p} + b^* b + c c^*) \gamma - \gamma  *( F_{A_q} + b b^* + c^* c ) \right) & = \lambda \gamma
\end{align*}
for some $\lambda < 0$. 

\end{enumerate}

\end{corollary}

Now specialize again to $\U(2,1)$, i.e.\ $\rank E_i=i$.

\medskip\noindent
{\bf (1)} 
$\mathcal{C}_a$. The negative eigenspace $\nu_a^-$ of the Hessian consists of holomorphic sections 
$$
\gamma \in 
H^0(E_1^*  E_2 \otimes K)
$$
 quiver bundle picture is
\begin{equation*}
(A)\qquad
\xygraph{
!{<0cm, 0cm>;<1.0cm, 0cm>:<0cm, 1.0cm>::}
!{(0,0) }*+{\bullet_{E_1}}="a"
!{(2,0) }*+{\bullet_{E_2}}="b"
"a" !{\ar @{-->}^{\gamma} "b"}
}
\end{equation*}

\medskip\noindent
{\bf (2)}  The critical points where $E_2 = S \oplus Q$ and the Higgs field is zero have negative eigenspace as follows.

\begin{enumerate}

\item $\mathcal{C}_{b_3}^{\ell}$. Since $\ell=d_S > d_1$, then the negative eigenspace $\nu^-_\ell$ of the Hessian consists of sections
\begin{equation*}
(\alpha, \beta_S, \gamma_Q) \in H^{0,1}(S^* Q) \oplus H^0(S^* E_1 \otimes K)\oplus   H^0( E_1^* Q \otimes K) .
\end{equation*}
These show up in the quiver bundle picture as dashed arrows in the diagram below
\begin{equation*}
(B_3)\qquad
\xygraph{
!{<0cm, 0cm>;<1.0cm, 0cm>:<0cm, 1.0cm>::}
!{(2,-1) }*+{\bullet_{S}}="a"
!{(0,0) }*+{\bullet_{E_1}}="b"
!{(2,1) }*+{\bullet_{Q}}="c"
"a" !{\ar @{-->}^{\beta_S} "b"}
"b" !{\ar @{-->}^{\gamma_Q} "c"}
"a" !{\ar @{-->}_{\alpha} "c"}}
\end{equation*}

\item $\mathcal{C}_{b_1}^{\ell}$. Since $d_2/2<\ell=d_S < d_1$, then the negative eigenspace of the Hessian is as follows
\begin{equation*}
\nu^-_\ell = \left\{ (\alpha, \gamma_S, \gamma_Q) \in H^{0,1}(S^* Q) \oplus H^0( E_1^* S \otimes K) \oplus \Omega^0(E_1^* Q \otimes K) \, : \, \bar{\partial}_A \gamma_Q + \alpha\gamma_S = 0 \right\}
\end{equation*}
The quiver bundle diagram is
\begin{equation*}
(B_1)\qquad
\xygraph{
!{<0cm, 0cm>;<1.0cm, 0cm>:<0cm, 1.0cm>::}
!{(2,-1) }*+{\bullet_{S}}="a"
!{(0,0) }*+{\bullet_{E_1}}="b"
!{(2,1) }*+{\bullet_{Q}}="c"
"b" !{\ar @{-->}_{\gamma_S} "a"}
"b" !{\ar @{-->}^{\gamma_Q} "c"}
"a" !{\ar @{-->}_{\alpha} "c"}}
\end{equation*}

\noindent
This is the case where the negative eigenspace of the Hessian is a singular space, and not a vector space.

As for the example of stable pairs (see for example \cite[Lemma 8.3.12]{WentworthWilkin11}), the idea is to further decompose the negative eigenspace of the Hessian. We have the equations
\begin{equation}\label{eqn:B2-negative}
\bar{\partial}_{A_1}^* \alpha  = 0 \ , \ 
\bar{\partial}_A \gamma_Q + \alpha\gamma_S  = 0 \ , \
\bar{\partial}_A \gamma_S  = 0 .
\end{equation}
Consider  the projection  from the solutions of 
\eqref{eqn:B2-negative} to the set $\{ \gamma_S \, : \, \bar{\partial}_A \gamma_S = 0 \}$. The remaining two equations are linear in $(\alpha, \gamma_Q)$, and therefore the fibers of this projection are vector spaces. The goal is to compute the dimension of these fibers (which will depend on $\gamma_S$). 

The case where $\gamma_S= 0$ is easy, since the equations decouple and the space of solutions is
$H^{0,1}(S^* Q) \oplus H^{1,0}(E_1^* Q)$.
When $\gamma_S \neq 0$ then the equations do not decouple, and we need to consider the following deformation complex
\begin{equation}\label{eqn:B2-def-complex}
\Omega^{0,1}(S^* Q) \oplus \Omega^{1,0}(E_1^* Q) \stackrel{D}{\xrightarrow{\hspace*{.75cm}}} \Omega^0(S^* Q) \oplus \Omega^{1,1}(E_1^* Q) ,
\end{equation}
where $D(\alpha, \gamma_Q) = (\bar{\partial}_A^* \alpha, \bar{\partial} \gamma_Q + \alpha\gamma_S)$. The adjoint  is
$
D^* (u, \eta) = \left( \bar{\partial}_A u - \bar{*} \left( \gamma_S \bar{*} \eta \right) , \bar{\partial}_A^* \eta \right) 
$.
We claim that the kernel of $D^*$ is trivial. To see this, note that $(u, \eta) \in \ker D^*$ implies that $\bar{\partial}_A^* \eta = 0$, and $\bar{\partial}_A u - \bar{*} \left( \gamma_S \bar{*} \eta \right) = 0$. The following calculation shows that this second equation decouples
\begin{align*}
\left< \bar{\partial}_A u, \bar{*} \left( \gamma_S \bar{*} \eta \right) \right> & = \left< u, \bar{\partial}_A^* \bar{*} \left( \gamma_S\bar{*} \eta \right) \right> 
  = \left< u, - \bar{*} \bar{\partial}_A \bar{*} \bar{*} \left( (\gamma_S) \bar{*} \eta \right) \right> \\
 & = \left< u, \bar{*} \bar{\partial}_A \left( \gamma_S \bar{*} \eta \right) \right> 
  = \left< u, \bar{*} \left( \bar{\partial}_A \gamma_S \bar{*} \eta - \gamma_S \bar{\partial}_A \bar{*} \eta \right) \right> = 0 
\end{align*}
since $\bar{\partial}_A \gamma_S = 0$ and $- \bar{*} \bar{\partial}_A \bar{*} \eta = \bar{\partial}_A^* \eta = 0$ by assumption.  

Therefore $(u, \eta) \in \ker D^*$ implies that
$
\bar{\partial}_A u  = 0
$,
$\bar{*} \left( \gamma_S \bar{*} \eta \right)  = 0
$, and
$
\bar{\partial}_A^* \eta  = 0
$.
Since $\deg S^* Q = 0$ then the first equation implies that $u=0$. Since $\gamma_S \neq 0$ then the second equation implies that $\eta = 0$, and together this shows that the kernel of $D^*$ is trivial. Therefore we can compute $\dim_\CBbb \ker D$ from the index of the complex and Riemann-Roch
\begin{align*}
\dim_\CBbb \ker D & = h^{0,1}(S^* Q) - h^0(S^* Q) + h^{1,0}(S^* Q) - h^{1,1}(E_1^* Q) \\
 & = g-1 - \deg (S^* Q) + h^{0,1}(Q^* E_1) - h^0(Q^* E_1) \\
 & = g-1 - \deg(S^* Q) + g-1-\deg(Q^* E_1) \\
 & = 2g-2 + d_S - d_1 
\end{align*}

We have proven the following
\begin{lemma}
The projection map  from the space of  solutions to
 \eqref{eqn:B2-negative} to the set 
$\{ \gamma_S \, : \, \bar{\partial}_A \gamma_S = 0 \}$ has linear fibers. 
The fiber over zero is isomorphic to $H^{0,1}(S^* Q) \oplus
H^{1,0}(E_1^* Q)$,
 and the fiber over any nonzero point has dimension $2g-2+d_S - d_1$.
\end{lemma}

\item When $\ell = d_1$ then the homomorphism $\gamma_S$ in the
diagram $(B_1)$ no longer corresponds to a negative eigenvalue of the Hessian (the eigenvalue is now zero). Therefore we have
\begin{equation*}
\nu^-_{\ell} = H^{0,1}(S^* Q) \oplus H^0(E_1^* Q \otimes K) 
\end{equation*}
and the quiver bundle picture is
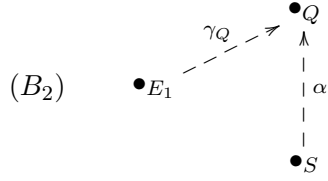
\begin{equation*}
(B_2)\qquad
\xygraph{
!{<0cm, 0cm>;<1.0cm, 0cm>:<0cm, 1.0cm>::}
!{(2,-1) }*+{\bullet_{S}}="a"
!{(0,0) }*+{\bullet_{E_1}}="b"
!{(2,1) }*+{\bullet_{Q}}="c"
"b" !{\ar @{-->}^{\gamma_Q} "c"}
"a" !{\ar @{-->}_{\alpha} "c"}}
\end{equation*}

\end{enumerate}

\medskip\noindent{\bf (3)} 
 The critical points where $E_2 = S \oplus Q$ and the Higgs field is nonzero have negative eigenspace as follows.

\begin{enumerate}

\item $\mathcal{C}_{c_1}^{\ell}$. 
 Since $\ell= d_S > \frac{1}{2}(d_Q + d_2)$ then the negative eigenspace $\eta^-_\ell$ of the Hessian is
\begin{equation*}
\eta^-_\ell = H^{0,1}(S^* Q) \oplus H^0(S^* E_1 \otimes K)  
\end{equation*}

The quiver bundle picture is
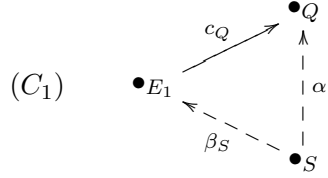
\begin{equation*}
(C_1)\qquad 
\xygraph{
!{<0cm, 0cm>;<1.0cm, 0cm>:<0cm, 1.0cm>::}
!{(2,-1) }*+{\bullet_{S}}="a"
!{(0,0) }*+{\bullet_{E_1}}="b"
!{(2,1) }*+{\bullet_{Q}}="c"
"a" !{\ar @{-->}^{\beta_S} "b"}
"a" !{\ar @{-->}_{\alpha} "c"}
"b" !{\ar @{>}^{c_Q} "c"}   }
\end{equation*}

\item $\mathcal{C}_{c_2}^{\ell}$. Since $d_Q < \frac{1}{2} ( d_S + d_1)$ then the negative eigenspace $\zeta^-_\ell$ of the Hessian consists of pairs $(\alpha, \gamma_Q) \in \Omega^{0,1} (S^* Q) \oplus \Omega^0(E_1^* Q \otimes K)$ such that
$$
\bar{\partial}_{A_1}^* \alpha - \bar{*} ( c_S (\bar{*} \gamma_Q))  =
0 \ ,\ 
\bar{\partial}_A \gamma_Q + \alpha c_S   =  0
$$
The quiver bundle picture is
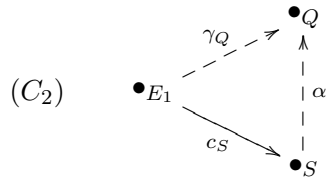
\begin{equation*}
(C_2)\qquad
\xygraph{
!{<0cm, 0cm>;<1.0cm, 0cm>:<0cm, 1.0cm>::}
!{(2,-1) }*+{\bullet_{S}}="a"
!{(0,0) }*+{\bullet_{E_1}}="b"
!{(2,1) }*+{\bullet_{Q}}="c"
"b" !{\ar @{-->}^{\gamma_Q} "c"}
"a" !{\ar @{-->}_{\alpha} "c"}
"b" !{\ar @{>}_{c_S} "a"    } }
\end{equation*}

\noindent
Note that these equations are both linear in $(\alpha, \gamma_Q)$, and that they correspond to the harmonic forms in the middle term of the following deformation complex
\begin{equation*}
\Omega^0(S^* Q) \stackrel{D_1}{\xrightarrow{\hspace*{.75cm}}}
 \Omega^{0,1}(S^* Q) \oplus \Omega^{1,0}(E_1^* Q)
\stackrel{D_2}{\xrightarrow{\hspace*{.75cm}}} \Omega^{1,1}(E_1^* Q) 
\end{equation*} 
where the maps $D_1$ and $D_2$ are
$$
D_1(u)  = \left( \bar{\partial}_A u, -u c_S \right) \ , \ 
D_2(\alpha, \gamma)  = \bar{\partial}_A \gamma_Q + \alpha c_S 
$$
(A calculation shows that $D_2 \circ D_1 = 0$.) The corresponding adjoints are
$$
D_1^* (\alpha, \gamma)  = \bar{\partial}_{A_1}^* \alpha - \bar{*} 
( c_S \bar{*} \gamma_Q) \ , \ 
D_2^* (\eta)  = \left( \bar{*}(c_S \bar{*} \eta), 
\bar{\partial}_A^* \eta \right) 
$$
If $c_S \neq 0$, then the maps $u \mapsto -uc_S$ 
and $\eta \mapsto \bar{*}(c_S\bar{*} \eta)$ both have trivial kernel, 
and hence the dimension of the harmonic forms 
in the middle term is equal to the index of the complex.
\begin{align*}
\dim_\CBbb \left(  \ker D_1^* \cap \ker D_2 \right) & = h^{0,1}(L_1^* L_2) - h^0(L_1^* L_2) + h^{1,0}(E_1^* L_2) - h^{1,1}(E_1^* L_2) \\
 & = g-1 + \ell - \ell_2 + h^{0,1}(L_2^* E_1) - h^0(L_2^* E_1) \\
 & = g-1+\ell - \ell_2 + g-1 + \ell_2 - d_2 \\
 & = 2g-2 + \ell - d_2 
\end{align*}
Therefore the negative eigenspace of the Hessian at these critical points has constant 
(complex) dimension $2g-2 + \ell - d_2$.

\item $\mathcal{C}_{c_3}^{\ell}$. Since $d_Q < \frac{1}{2} (d_S + d_2)$ then the negative eigenspace $\zeta^-_\ell$ of the Hessian consists of pairs $(\alpha, \gamma_Q) \in \Omega^{0,1} (S^* Q) \oplus \Omega^0(E_1^* Q \otimes K)$ such that
$
\bar{\partial}_{A_1}^* \alpha  = 0
$,
$
\bar{\partial}_A \gamma_Q  = 0 
$,
and so the space of solutions is isomorphic to $H^{0,1}(S^* Q) \oplus H^0(E_1^* Q \otimes K)$. The quiver bundle picture is
\begin{equation*}
(C_3)\qquad
\xygraph{
!{<0cm, 0cm>;<1.0cm, 0cm>:<0cm, 1.0cm>::}
!{(2,-1) }*+{\bullet_{S}}="a"
!{(0,0) }*+{\bullet_{E_1}}="b"
!{(2,1) }*+{\bullet_{Q}}="c"
"b" !{\ar @{-->}^{\gamma_Q} "c"}
"a" !{\ar @{-->}_{\alpha} "c"}
"a" !{\ar @{>}^{b_S} "b"} }
\end{equation*}

\end{enumerate}

\subsection{Cohomology of negative normal directions}\label{sec:neg-directions} 
We now compute the relative cohomology groups of the negative normal
spaces given in the previous section.  

\medskip\noindent {\bf (1)}   Consider the case of the A-stratum, where $\ell=d_2/2$. Let
\begin{align*}
\nu_a^-&=\left\{ (E_1, E_2, 0, \gamma) : \gamma\in H^0(E_1^\ast E_2 \otimes K)\ , \ \text{$E_2$ semistable}\right\} \\
\nu_a'&=\left\{ (E_1, E_2, 0, \gamma)\in \nu_a^- : \gamma\neq 0\right\}
\end{align*}
We also fix
\begin{equation} \label{eqn:sigma-min}
\sigma_{min} := 2g-2-d_1+d_2/2+ 1/4
\end{equation}
The important point is that $ \tfrac{1}{2}\deg(E_1^\ast E_2 \otimes K)<\sigma_{min}< \lfloor \tfrac{1}{2}\deg(E_1^\ast E_2 \otimes K)\rfloor+1$.

\begin{lemma} \label{lem:neg-a-cohomology}
For the A stratum:
\begin{enumerate}
\item $H^\ast_\G(\nu_a^-)\simeq H^\ast_\G(\Acal(E_1)\times \Acal^{ss}(E_2))$
\item  If $d_2$ is odd, then
$H^\ast_\G(\nu_a')\simeq H^\ast(\Ncal_{\sigma_{min}}(E_1^\ast E_2\otimes K))$
\item If $d_2$ is even, then
$$
H^\ast_\G(\nu_a')\simeq H^\ast(\Ncal_{\sigma_{min}}(E_1^\ast E_2\otimes K))
\oplus H^{\ast-2(2g-2-d_1+d_2/2)}_{S^1\times S^1}(J(X)\times J(X)\times S^{2g-2-d_1+d_2/2}X)
$$
\end{enumerate}
\end{lemma}
\begin{proof}
Part (i) follows by the deformation retraction $\gamma\mapsto 0$.  
Let $E=E_1^\ast E_2\otimes K$.
Part (ii) follows because $\Gcal(E)$ acts freely with 
 $\nu_a'/\G=\Ncal_{\sigma_{min}}(E)$.  
 Part (iii) is slightly more subtle.  A $\sigma_{min}$-Bradlow stable pair is a nonvanishing section 
 $\gamma\in H^0(E)$ with the additional assumption, in case $E$ is strictly semistable, 
that $\gamma$ does not lie in the maximally destabilizing subbundle. 
 Hence, the space $\nu_a'$ is obtained  by attaching the first nonminimal stratum to the Bradlow semistable stratum in the $\sigma_{min}$-YMH stratification of the space of pairs given in \cite[Section 8.2.1]{WentworthWilkin11}.  Then part (iii) follows from the computation in 
 \cite[Theorem 8.4.1]{WentworthWilkin11}. 
\end{proof}

\medskip\noindent {\bf (2)} 
 $\tfrac{1}{3}(d_1+d_2)< \ell \leq d_2-d_1+2g-2$.
This is the case of the  $C_1$-stratum; see the quiver diagram
$(C_1)$.  Define the following spaces
\begin{align*}
\eta^-_{\ell} &=\left\{ (\alpha,\beta_S) :  \dbar^\ast\alpha=0\ ,\ \dbar\beta_S=0
\right\} \\
\eta'_{\ell} &=\left\{ (\alpha,\beta_S) \in \eta^-_{\ell}: (\alpha,\beta_S)\neq 0\right\} \\
\eta''_{\ell} &=\left\{ (\alpha,\beta_S) \in \eta^-_{\ell}: \alpha\neq 0\right\}  
\end{align*}
Then by the argument in \cite{DWWW11} we have
\begin{lemma}
For the $C_1$ stratum,
\begin{align}
H^\ast_\G(\eta^-_{\ell}, \eta''_{\ell})
 &= H^{\ast-2(d_2-2\ell+g-1)}_{S^1\times S^1}
\left(   J(X)\times J(X)\times S^{\ell-d_1+2g-2}X          \right)                 \label{eqn:zeta''_C1} \\
H^\ast_\G(\eta'_{\ell},\eta''_{\ell} )
&=H^{\ast-2(d_2-2\ell+g-1)}_{S^1}\left(   J(X)\times S^{d_2-d_1+2g-2-\ell}X  \times
S^{d_1-\ell+2g-2}X    \right) \label{eqn:zeta'_C1}  \\
H^\ast_\G(\eta^-_{\ell}, \eta''_{\ell}) &=
H^\ast_\G(\eta^-_{\ell}, \eta'_{\ell})   \oplus
H^\ast_\G(\eta'_{\ell},\eta''_{\ell} )            \label{eqn:zeta_C1}   
\end{align}
\end{lemma}

\medskip\noindent {\bf (3)}
$\tfrac{1}{3}(2d_2-d_1)<\ell<d_1$. These are the $B_1$ and $C_2$ strata.  
Consider first the diagram $(B_1)$.
Define the following spaces
\begin{align*}
\nu^-_{\ell} &=\left\{ (\alpha,\gamma_S,\gamma_Q) :  \dbar^\ast\alpha=0\ ,\ \dbar\gamma_Q+\alpha\gamma_S=0\ ,\ \dbar\gamma_S=0
\right\} \\
\nu'_{\ell} &=\left\{ (\alpha,\gamma_S,\gamma_Q) \in \nu^-_{\ell}: (\alpha,\gamma_S,\gamma_Q)\neq 0\right\} \\
\nu''_{\ell} &=\left\{ (\alpha,\gamma_S,\gamma_Q) \in \nu^-_{\ell}: \alpha\neq 0\right\}  \\
\omega_\ell & = \left\{  (\alpha,\gamma_S,\gamma_Q) \in \nu^-_{\ell} \, : \, (\alpha, \gamma_Q)\neq 0 \right\}  
\end{align*}
\begin{lemma} \label{lem:b1}
For the $B_1$ stratum, $\ell\leq d_2-d_1+2g-2$,
\begin{align}
H^\ast_\G(\nu^-_\ell, \nu''_\ell) &=
 H^{\ast-2(2\ell-d_2+g-1)}_{S^1\times S^1\times S^1}\left(   J(X)\times J(X)\times J(X)               \right)                 \label{eqn:nu''_B1} \\
H^\ast_\G(\nu'_\ell, \nu''_\ell) &=   H^\ast_\G(\nu'_\ell, \omega_\ell)   \oplus      H^\ast_\G(\omega_\ell,\nu''_\ell )            \label{eqn:nu'_B1}   \\
H^\ast_\G(\nu'_\ell, \omega_\ell) &=
H^{\ast-2(\ell-d_1+2g-2)}_{S^1\times S^1}\left(   J(X)\times J(X)\times
 S^{\ell-d_1+2g-2}X    \right)
\label{eqn:omega'_B1} \\
H^\ast_\G(\omega_\ell,\nu''_\ell ) &=
H^{\ast-2(2\ell-d_2+g-1)}_{S^1\times S^1}\left(   J(X)\times J(X)\times
 S^{d_2-d_1+2g-2-\ell}X    \right)
\label{eqn:omega''_B1} 
\end{align}
If $d_2-d_1+2g-2<\ell<d_1$, then \eqref{eqn:nu''_B1} holds, with
\begin{align}
H^\ast_\G(\nu^-_\ell, \nu''_\ell) &=   H^\ast_\G(\nu^-_\ell, \nu'_\ell)   \oplus
H^\ast_\G(\nu'_\ell,\nu''_\ell )            \label{eqn:nu'-B1-large}   \\
H^\ast_\G(\nu'_\ell, \nu''_\ell) &=
H^{\ast-2(\ell-d_1+2g-2)}_{S^1\times S^1}\left(   J(X)\times J(X)\times
 S^{\ell-d_1+2g-2}X    \right)
\label{eqn:nu''-B1-large}
\end{align}
\end{lemma}

\begin{proof}
Notice that \eqref{eqn:nu''_B1} follows by retracting $(\gamma_S,\gamma_Q)\mapsto 0$ and using the Atiyah-Bott argument. 
 Consider the following commutative diagram.

\begin{equation}\label{eqn:type-B-diagram}
\begin{split}
\xymatrix{
&& \vdots \ar[d] &&&\\
&\cdots \ar[r] & H_\mathcal{G}^p(\nu^-_{\ell}, \nu'_{\ell}) \ar[r] \ar[d] & H_\mathcal{G}^p(\nu^-_{\ell}) \ar[r] & H_\mathcal{G}^p(\nu'_{\ell}) \ar[r] & \cdots \\
 && H_\mathcal{G}^p(\nu^-_{\ell}, \nu''_{\ell}) \ar[d] \ar[ur]^\xi \ar[dr]^{\xi''} &&&\\
\cdots\ar[r] & H^p_\G(\nu'_{\ell},\omega_\ell) \ar[r] & H_\mathcal{G}^p(\nu'_{\ell}, \nu''_{\ell}) \ar[d] \ar[r]^\beta & H^p_\G(\omega_\ell, 
\nu''_{\ell}) \ar[r] & \cdots &
\\
 && \vdots &&&
}
\end{split}
\end{equation}
By the argument in \cite{DWWW11} and assuming \eqref{eqn:omega''_B1}, 
the map $\xi''$ is surjective. It follows that the lower 
horizontal exact sequence splits.  Thus, 
\eqref{eqn:nu'_B1} follows from \eqref{eqn:omega''_B1}.  Define 
the following spaces
$$
\negthickspace
\begin{array} {lclclcl}
W_\ell &=& \left\{  (\alpha,\gamma_S,\gamma_Q) \in \nu^-_{\ell}  
:  \gamma_S = 0 \right\} && 
W_\ell' & = & \left\{  (\alpha,\gamma_S,\gamma_Q) \in \nu_{\ell} : 
 \gamma_S = 0, (\alpha, \gamma_Q) \neq 0 \right\} \\
W_\ell'' & = &\left\{  (\alpha,\gamma_S,\gamma_Q) \in
\nu^-_{\ell} :  \gamma_S = 0, \alpha \neq 0 \right\} && 
Z_\ell & = & \left\{  (\alpha,\gamma_S,\gamma_Q) \in \nu^-_{\ell} :  (\gamma_S,\gamma_Q) = 0, \alpha\neq 0 \right\} \\
R_\ell &=& \left\{  (\alpha,\gamma_S,\gamma_Q) \in \nu^-_{\ell} :
 \gamma_Q \neq 0, (\alpha, \gamma_S) = 0 \right\} && 
Y_\ell'' & = &\left\{  (\alpha,\gamma_S,\gamma_Q) \in \nu^-_{\ell}  :  \gamma_S \neq 0, (\alpha, \gamma_Q) \neq 0 \right\} \\
T_\ell &=& \left\{  (\alpha,\gamma_S,\gamma_Q) \in \nu^-_{\ell} :
 \gamma_S \neq 0, (\alpha, \gamma_Q) = 0 \right\} && 
Y_\ell' & = & \left\{  (\alpha,\gamma_S,\gamma_Q) \in \nu^-_{\ell} :  \gamma_S \neq 0 \right\} \\
Y_\ell'' & = & \left\{  (\alpha,\gamma_S,\gamma_Q) \in \nu^-_{\ell}
 :  \gamma_S \neq 0, (\alpha, \gamma_Q) \neq 0 \right\} && 
T_\ell &=& \left\{  (\alpha,\gamma_S,\gamma_Q) \in \nu^-_{\ell} :  \gamma_S \neq 0, (\alpha, \gamma_Q) = 0 \right\} 
\end{array}
$$
\noindent 
Note that $Y_\ell' = \nu^-_{\ell} \setminus W_\ell = \nu'_{\ell} \setminus 
W_\ell'$ and $Y_\ell'' = \omega_{\ell} \setminus W_\ell'$.
By the retraction $\gamma_S\mapsto 0$, the pair $(\omega_\ell, \nu''_\ell)\simeq (W'_\ell, W''_\ell)$.
By excision,
$$
H^\ast_\G(W'_\ell, W''_\ell)\simeq H^\ast_\G(W'_\ell\setminus Z_\ell, W''_\ell\setminus Z_\ell)
$$
Now $W'_\ell\setminus Z_\ell$ fibers over $R_\ell$ with fiber 
dimension $d_S-d_Q+g-1$.  Hence, \eqref{eqn:omega''_B1} follows from the 
Thom isomorphism.  Finally, for \eqref{eqn:omega'_B1}
we need the following lemma, whose proof is straightforward. 

\begin{lemma} \label{lem:gamma_S}
For fixed $\gamma_S\neq 0$, the space of solutions $(\alpha, \gamma_Q)$ to
$
 \dbar^\ast\alpha=0$,
 $
 \dbar\gamma_Q+\alpha\gamma_S=0
 $,
 has dimension $= \ell-d_1+2g-2$.
\end{lemma}

\noindent
Excision of $W'_\ell$ gives
$
H^\ast_\G(\nu'_\ell, \omega_\ell)\simeq H^\ast_\G(\nu'_\ell\setminus W'_\ell, \omega_\ell\setminus W'_\ell)=H^\ast_{\G}(Y'_\ell, Y''_\ell)
$.
Now by the lemma, $Y'_\ell$ fibers over $T_\ell$ with fiber dimension 
$\ell-d_1+2g-2$, and \eqref{eqn:omega'_B1} again follows from Thom isomorphism.
In case $d_2-d_1+2g-2<\ell<d_1$, 
then notice that $W'_\ell$ is closed in $\nu''_\ell$.  Hence, \eqref{eqn:nu''-B1-large} follows by Lemma \ref{lem:gamma_S} and excision. Eq.\ \eqref{eqn:nu'-B1-large} follows by the argument in \cite{DWWW11}.
\end{proof}

For $C_2$, the normal directions are given by
\begin{align*}
\zeta^-_{\ell} &=\left\{ (\alpha,\gamma_Q) :  \dbar^\ast\alpha=0\ ,\ \dbar\gamma_Q+\alpha c_S=0
\right\} \\
\zeta'_{\ell} &=\left\{ (\alpha,\gamma_Q) \in \zeta^-_{\ell}: (\alpha,\gamma_Q)\neq 0\right\} 
\end{align*}
where $c_S\neq 0$. It follows from Lemma \ref{lem:gamma_S} that
\begin{equation} \label{eqn:c2}
H^\ast_\G(\zeta^-_\ell, \zeta'_\ell)\simeq
 H^{\ast-2(\ell-d_1+2g-2)}_{S^1\times S^1}(J(X)\times J(X)\times
S^{\ell-d_1+2g-2}X)
\end{equation}
\begin{remark} \label{rem:b1-c2}
For the $B_1$ and $C_2$ strata, $\ell\leq d_2-d_1+2g-2$, $H^\ast_\G(\nu'_\ell, \omega_\ell)\simeq H^\ast_\G(\zeta^-_\ell, \zeta'_\ell)$. In case $d_2-d_1+2g-2<\ell<d_1$, then
$H^\ast_\G(\nu'_\ell, \nu''_\ell)\simeq H^\ast_\G(\zeta^-_\ell, \zeta'_\ell)$.
\end{remark}

\medskip\noindent {\bf (4)} 
$d_1\leq \ell$.  These are the  $B_2, B_3$, and $C_3$ strata. 
Consider first the  the $(C_3)$ diagram.
Define the following spaces
\begin{align*}
\zeta^-_{\ell} &=\left\{ (\alpha,\gamma_Q) :  \dbar^\ast\alpha=0\ ,\ \dbar\gamma_Q=0
\right\} \\
\zeta'_{\ell} &=\left\{ (\alpha,\gamma_Q) \in \zeta^-_{\ell}: (\alpha,\gamma_Q)\neq 0\right\} \\
\zeta''_{\ell} &=\left\{ (\alpha,\gamma_Q) \in \zeta^-_{\ell}: \alpha\neq 0\right\}  
\end{align*}
Then by the argument in \cite{DWWW11} we have
\begin{lemma} \label{lem:C_3}
For the $C_3$ stratum, if $\ell\leq d_2-d_1+2g-2$ then
\begin{align}
H^\ast_\G(\zeta^-_\ell, \zeta''_\ell) &=
H^{\ast-2(2\ell-d_2+g-1)}_{S^1\times S^1}\left(   J(X)\times
J(X)\times S^{d_1-\ell+2g-2}X          \right)                 \label{eqn:zeta''_C3} \\
H^\ast_\G(\zeta'_\ell,\zeta''_\ell ) &=H^{\ast-2(2\ell-d_2+g-1)}_{S^1}\left(
J(X)\times S^{d_2-d_1-\ell+2g-2}X  \times S^{d_1-\ell+2g-2}X    \right)
\label{eqn:zeta'_C3}  \\
H^\ast_\G(\zeta^-_\ell, \zeta''_\ell) &=   H^\ast_\G(\zeta^-_\ell, \zeta'_\ell)   \oplus      H^\ast_\G(\zeta'_\ell,\zeta''_\ell )            \label{eqn:zeta_C3}   
\end{align}
If $d_2-d_1+2g-2<\ell\leq d_1+2g-2$ then 
\begin{equation}  \label{eqn:c3-bundle-case}
H^\ast_\G(\zeta^-_\ell, \zeta'_\ell)\simeq
 H^{\ast-2(2\ell-d_2+2g-2)}_{S^1\times S^1}(J(X)\times J(X)\times
S^{d_1-\ell+2g-2}X)
\end{equation}
\end{lemma}

The $B_2$ case is exactly the same as the $C_3$ case. We define the spaces $\nu^-_{d_1}$,
$\nu'_{d_1}$, and 
$\nu''_{d_1}$ by analogy to $\zeta^-_\ell$, $\zeta'_\ell$, and $\zeta''_\ell$ above.
\begin{lemma} \label{lem:b2}
For the $B_2$ stratum, if $d_1\leq d_2-d_1+2g-2$, then
\begin{align}
H^\ast_\G(\nu^-_{d_1}, \nu''_{d_1}) &=
 H^{\ast-2(2d_1-d_2+g-1)}_{S^1\times S^1}\left(   J(X)\times J(X)\times J(X)  \right)
\label{eqn:nu''_B2} \\
H^\ast_\G(\nu'_{d_1},\nu''_{d_1} ) &=
H^{\ast-2(2d_1-d_2+g-1)}_{S^1}\left(   J(X)\times J(X)\times S^{d_2-2d_1+2g-2}X  \right) 
\label{eqn:nu'_B2}  \\
H^\ast_\G(\nu^-_{d_1}, \nu''_{d_1}) &=  
 H^\ast_\G(\nu^-_{d_1}, \nu'_{d_1})   \oplus   
   H^\ast_\G(\nu'_{d_1},\nu''_{d_1} )            \label{eqn:nu_B2}   
\end{align}
If $d_2-d_1+2g-2< d_1$, then
\begin{equation} \label{eqn:b2-large-degree}
H^\ast_\G(\nu^-_{d_1}, \nu'_{d_1}) \simeq H^{\ast-2(2d_1-d_2+d-1)}_{S^1\times S^1\times S^1}(J(X)\times J(X)\times J(X))
\end{equation}
\end{lemma}

Finally,  consider  the  $(B_3)$ diagram.  There are three cases.  First, if $d_1<\ell\leq d_2-d_1+2g-2$,
define the following spaces
\begin{align*}
\nu^-_{\ell} &=\left\{ (\alpha,\beta_S,\gamma_Q) :  \dbar^\ast\alpha=0\ ,\ \dbar\beta_S=0\ ,\ \dbar\gamma_Q=0
\right\} \\
\nu'_{\ell} &=\left\{ (\alpha,\beta_S,\gamma_Q) \in \nu^-_{\ell}: (\alpha,\beta_S,\gamma_Q)\neq 0\right\} \\
\nu''_{\ell} &=\left\{ (\alpha,\beta_S,\gamma_Q) \in \nu^-_{\ell}: \alpha\neq 0\right\}  \\
\omega_\ell & = \left\{  (\alpha,\beta_S,\gamma_Q) \in \nu^-_{\ell} \, : \, (a, \gamma_Q)\neq 0 \right\}  
\end{align*}
\begin{lemma} \label{lem:b3-case1}
For the $B_3$ stratum, $d_1< \ell\leq d_2-d_1+2g-2$,
\begin{align}
H^\ast_\G(\nu^-_\ell, \nu''_\ell) &=
 H^{\ast-2(2\ell-d_2+g-1)}_{S^1\times S^1\times S^1}\left(   J(X)\times J(X)\times J(X)               \right)                 \label{eqn:nu''_B3} \\
H^\ast_\G(\nu'_\ell, \nu''_\ell) &=   H^\ast_\G(\nu'_\ell, \omega_\ell)   \oplus      H^\ast_\G(\omega_\ell,\nu''_\ell )            \label{eqn:nu'_B3}   \\
H^\ast_\G(\nu'_\ell, \omega_\ell) &=H^\ast_\G(\zeta^-_\ell, \zeta_\ell') 
\label{eqn:omega'_B3} \\
H^\ast_\G(\omega_\ell,\nu''_\ell ) &=
H^{\ast-2(2\ell-d_2+g-1)}_{S^1\times S^1}\left(   J(X)\times
J(X)\times S^{d_2-d_1-\ell+2g-2}X    \right)
\label{eqn:omega''_B3} 
\end{align}
\end{lemma}
\begin{proof}  \eqref{eqn:nu''_B3} follows as before, and \eqref{eqn:nu'_B3} follows from \eqref{eqn:omega''_B3}.  For \eqref{eqn:omega'_B3},  use excision on the set $\{ \beta_S=0\}$.  Finally, for \eqref{eqn:omega''_B3}, first retract $\beta_S\mapsto 0$ and then excise the set $\{\gamma_S=0\}$.  The rest fibers over $\{\gamma_S\neq 0\}$, and the result follows from the Thom isomorphism.
\end{proof}

In case $d_2-d_1+2g-2<\ell\leq d_1+2g-2$, then $\gamma_Q\equiv 0$. 
Eq.\ \eqref{eqn:nu''_B3} holds as before, but now 
\begin{align}
\begin{split}\label{eqn:b3}
H^\ast_\G(\nu'_\ell, \nu''_\ell) &=  H^\ast_\G(\zeta^-_\ell, \zeta_\ell') \\
H^\ast_\G(\nu^-_\ell, \nu''_\ell) &=   H^\ast_\G(\nu^-_\ell, \nu'_\ell)   \oplus      H^\ast_\G(\nu'_\ell,\nu''_\ell )
\end{split}
\end{align}
  If $d_1+2g-2<\ell$, then both $\beta,\gamma\equiv 0$, and by the Atiyah-Bott isomorphism
\begin{equation} \label{eqn:b3-large-degree}
H^\ast_\G(\nu^-_\ell, \nu'_\ell) \simeq H^{\ast-2(2\ell-d_2+g-1)}_{S^1\times S^1\times S^1}(J(X)\times J(X)\times J(X))
\end{equation}

\subsection{The Morse-Bott Lemma} 
The goal of this section is prove the validity of the Morse-Bott isomorphism, which relates the equivariant cohomology of the pair of successive strata to the equivariant cohomology of the pair consisting of negative normal directions and nonzero negative normal directions.  Because of singularities, Bott's argument 
in \cite{Bott54} does not apply, and as in \cite{DWWW11} and \cite{WentworthWilkin11} we need to 
circumvent this. 
 In fact, we do not prove the Morse-Bott lemma for all critical sets.  Nevertheless, the results below are sufficient for the cohomological calculations in the next section.

We begin with particular regions of the parameter $\ell\in \Delta_{d_1,d_2}$ using the definition on page \pageref{def:regions}.

\begin{proposition} \label{prop:bott1}
 For regions $({\bf II})$ and $({\bf III})$, 
\begin{align}
H^\ast_\G(X_{d_2/2}^\ast\cup \Scal_a, X_{d_2/2}^\ast) &\simeq  H^\ast_\G(\nu^-_a, \nu_a') 
\label{eqn:mb-a} \\
H^\ast_\G(X_\ell, X_\ell') &\simeq  H^\ast_\G(\nu^-_\ell, \nu_\ell') 
\label{eqn:mb-b} \\
H^\ast_\G(X_\ell', X^\ast_{\ell}) &\simeq  H^\ast_\G(\zeta^-_\ell, \zeta_\ell')   \qquad (\ell\leq d_1+2g-2)
\label{eqn:mb-c}
\end{align}
Moreover, in these regions the inclusions $X_\ell'\subset X_\ell$ and $X_\ell^\ast\subset X_\ell'$
 induce surjections 
$
H^\ast_\G(X_\ell)\lra H^\ast_\G(X_\ell')
$
and $
H^\ast_\G(X_\ell')\lra H^\ast_\G(X_\ell^\ast)
$.
\end{proposition}

We will need the following result.
Consider $\U(2,1)$ bundles where $b\equiv 0$, i.e.\ quiver bundles of the form
\begin{equation} \label{eqn:b-zero-quiver}
\xygraph{
!{<0cm, 0cm>;<1.0cm, 0cm>:<0cm, 1.0cm>::}
!{(0,0) }*+{\bullet_{E_1}}="a"
!{(2,0) }*+{\bullet_{E_2}}="b"
"a" !{\ar @{->}^{c} "b"}}
\end{equation}
The data is clearly equivalent to a choice of holomorphic section (also denoted $c$) of the bundle 
$E_1^\ast E_2\otimes K$.  We have the following
\begin{lemma} \label{lem:bradlow-stability}
For quivers of the type above, Higgs (semi)stability of $(E_2\oplus E_1, 0,c)$ is
 equivalent to Bradlow (semi)stability of the pair $(E_1^\ast E_2\otimes K, c)$
  for $\sigma=\sigma(d_1,d_2)$ as defined in \eqref{eqn:sigma}.
\end{lemma}

\begin{proof}
Set $\displaystyle\Phi=\left(\begin{matrix} 0&c\\ 0&0\end{matrix}\right)$ and $E=E_1^\ast E_2\otimes K$.
  Any line subbundle $S\subset E_2$ is automatically $\Phi$-invariant, so  Higgs semistability implies
$
d_S\leq \tfrac{1}{3}(d_1+d_2)
$.
If moreover $c(E_1)\subset S\otimes K$, then Higgs semistability implies 
$
\tfrac{1}{2}(d_S+d_1)\leq \tfrac{1}{3}(d_1+d_2)
$.
  On the other hand, $S\subset E_2$ gives a line subbundle $E_1^\ast S\otimes K\subset E$.  Then $\sigma$-semistability implies 
$
\deg(E_1^\ast S\otimes K) \leq \sigma$, or
$
d_S-d_1+ 2g-2 \leq \sigma
$.
If $c(E_1)\subset S\otimes K$, then the corresponding section of $E$ lies in 
$E_1^\ast S\otimes K\subset E$, so $\sigma$-semistability implies
$$
\sigma \leq \deg E -\deg(E_1^\ast S\otimes K\subset E) \leq d_2-d_1-d_S+2g-2
$$
Now for the given choice $\sigma=\sigma(d_1,d_2)$ as in \eqref{eqn:sigma}, the conditions for Higgs and $\sigma$-semistability are  equivalent.
\end{proof}

\begin{proof}[Proof of Proposition \ref{prop:bott1}]
For the $C_2$ stratum,  the $C_3$ stratum
$d_2-d_1+2g-2<\ell<d_1+2g-2$, 
 the $B_3$ stratum $d_1+2g-2<\ell$, or the $B_2$ stratum when $d_2-d_1+2g-2<d_1$,
 the negative normal directions are vector bundles.  The result then follows from a standard argument.  
  To prove \eqref{eqn:mb-b} for the portion of the  $B_1$ stratum where $d_2/2<\ell\leq \tfrac{1}{3}(d_1+d_2)$ (or $d_2-d_1+2g-2<\ell<d_1$),
 define the map $\pr : \Bcal(d_1,d_2)\to \Acal(E_2)$ by projection to the holomorphic structure on $E_2$.
 Let
 $$
K_\ell=\bigcup_{j>\ell} X_\ell\cap \pr^{-1}(\A_j(E_2))
$$
Then $K_\ell\subset X'_\ell$ is closed in $X_\ell$.  Hence, by excision,
$$
H^\ast_\Gcal(X_\ell, X'_\ell)\simeq H^\ast_\Gcal(X_\ell\setminus K_\ell, X'_\ell\setminus K_\ell)
$$
Moreover, the pair $(X_\ell\setminus K_\ell, X'_\ell\setminus K_\ell)$ is invariant under the 
scaling $b\mapsto 0$.  The same is true of the pair $(\nu_\ell^-, \nu'_\ell)$.  
Eq.'s \eqref{eqn:mb-a} and \eqref{eqn:mb-b} therefore reduce to the corresponding result 
for pairs \eqref{eqn:b-zero-quiver},  and hence they follow from 
 Lemma \ref{lem:bradlow-stability} and   \cite[eq.\ (8.28) and Sect.\ 8.3.6]{WentworthWilkin11}.

It remains to prove \eqref{eqn:mb-b} for the portion of the  $B_3$ stratum where $
\max\{d_1, d_2-d_1+2g-2\}<\ell\leq d_1+2g-2$.  For all integers $\ell>d_2/2$, 
let $X_\ell''=X_\ell\setminus\pr^{-1}(\Acal_\ell(E_2))$.
  Then it follows as in \cite[eq.\ (21)]{DWWW11} that
\begin{equation} \label{eqn:ab}
H^\ast_\Gcal(X_\ell, X''_\ell)\simeq H^\ast_\Gcal(\nu_\ell^-, \nu''_\ell)
\end{equation}
and by the Atiyah-Bott lemma, $H^\ast_\Gcal(X_\ell, X''_\ell)\to H^\ast_\Gcal(X_\ell)$ is injective.
We claim that for $k>d_2-d_1+2g-2$, $X_\ell''=X_\ell^\ast$.  
Indeed, it suffices to show that if $(E_2\oplus E_1,b,c)$ is semistable, then the Harder-Narasimhan type 
of $E_2$ is at most $d_2-d_1+2g-2$.  Suppose not and let $0\to S\to E_2\to Q\to 0$ be 
the Harder-Narasimhan filtration with $\deg S=\ell$.  
Then if $\ell>d_2-d_1+2g-2$, the induced map $c:E_1\to Q$ vanishes and $S\oplus E_1$ is $\Phi$-invariant. 
Hence,
$$
\tfrac{1}{2}(\ell+d_1)\leq \tfrac{1}{3}(d_1+d_2)\ \Longrightarrow\
\tfrac{1}{2}(d_2+2g-2)< \tfrac{1}{3}(d_1+d_2)\ \Longrightarrow\
2g-2<\tfrac{1}{3}(2d_1-d_2)\leq g-1
$$
where the last inequality comes from the bound on the Toledo invariant.  
This contradicts the assumption on the genus, and the claim follows.  Now 
the proof of \eqref{eqn:mb-b} follows from the fact that $H^\ast_\Gcal(X_\ell', X_\ell^\ast)\simeq 
H^\ast_\Gcal(\nu_\ell', \nu_\ell'')$ by \eqref{eqn:b3}, and the Five Lemma 
applied to the long exact sequence of the triple $(X_\ell, X_\ell', X_\ell^\ast)$.
\end{proof}

\begin{corollary} \label{cor:b1}
For the $B_1$ stratum in the portion of region $({\bf II})$ where
$d_2/2<\ell\leq \tfrac{1}{3}(d_1+d_2)$,
$
H^\ast_\Gcal(X_\ell, X^\ast_\ell)\simeq \ker \xi''
$.
If $d_2-d_1+2g-2<\ell<d_1$, $
H^\ast_\Gcal(X_\ell, X^\ast_\ell)\simeq H^\ast_\Gcal(\nu^-_\ell, \nu''_\ell)
$.
\end{corollary}

\begin{proof}
By the exact sequence of the triple $(X_\ell, X'_\ell, X^\ast_\ell)$, Remark \ref{rem:b1-c2}, and Proposition 
\ref{prop:bott1},
$$
\xymatrix{
&\cdots \ar[r] & H_\mathcal{G}^\ast(X_\ell, X_\ell') \ar[r] \ar[d]^{\begin{sideways}$\sim$\end{sideways}} &\ar[d]
 H_\mathcal{G}^\ast(X_\ell, X_{\ell}^\ast) \ar[r]&  \ar[d]^{\begin{sideways}$\sim$\end{sideways}}
H_\mathcal{G}^\ast(X'_\ell, X_{\ell}^\ast) \ar[r] & \cdots \\
&\cdots \ar[r] & H_\mathcal{G}^\ast(\nu^-_{\ell}, \nu'_{\ell}) \ar[r] & \ker\xi'' \ar[r]&
H_\mathcal{G}^\ast(\nu'_{\ell}, \omega_{\ell}) \ar[r] & \cdots 
}
$$
The first statement follows from the Five Lemma.  The proof of the second statement is similar.
\end{proof}

Now consider the region (${\bf I}$), 
which involves the  $C_1$ stratum.
We have the following

\begin{lemma}  \label{lem:C-bott}
For all $\tfrac{1}{3}(d_1+d_2)<\ell\leq  d_2-d_1+2g-2$,
 $H^\ast_\G(X^\ast_{\ell}, X''_\ell)\simeq H^\ast_\G(\eta'_{\ell}, \eta''_{\ell})$.
\end{lemma}

\begin{proof}
The argument is similar to the one in \cite[Section 3.1]{DWWW11}. 
Note that
 the set 
$
(X^\ast_\ell\setminus \Bcal^{ss}(d_1,d_2))\subset X''_\ell
$
is closed in $X^\ast_\ell$.  Hence, by excision
$$
H^\ast_\G(X^\ast_\ell, X''_\ell)\simeq H^\ast_\G( \Bcal^{ss}(d_1,d_2),
  \Bcal^{ss}(d_1,d_2)\setminus \pr^{-1}(\A_\ell(E_2)))
$$
By \cite{Wilkin08}, the YMH flow defines a $\G$-equivariant deformation retraction of the pair 
$$( \Bcal^{ss}(d_1,d_2),  \Bcal^{ss}(d_1,d_2)\setminus \pr^{-1}(\A_\ell(E_2)))$$
 with
$( \Bcal_{min}(d_1,d_2),  \Bcal_{min}(d_1,d_2)\setminus \pr^{-1}(\A_\ell(E_2)))$.
Note that $\Bcal_{min}(d_1,d_2)\cap \pr^{-1}(\A_\ell(E_2))$ lies in the smooth locus on which $\Gcal$ acts freely.
Excision then  reduces the computation  to Gothen's calculation in \cite{Gothen02}.
\end{proof}

By \eqref{eqn:ab},  Lemma \ref{lem:C-bott}, and  \eqref{eqn:zeta'_C1}, and the argument in \cite{DWWW11}, we have

\begin{corollary} \label{cor:C-bott}
For all $\tfrac{1}{3}(d_1+d_2)<\ell\leq d_2-d_1+2g-2$,
the map $H^\ast_\G(X_\ell, X_\ell'')\to H^\ast_\G(X^\ast_\ell, X_\ell'')$ is surjective.
\end{corollary}

\subsection{Proof of Theorem \ref{thm:perfection}}
\begin{lemma} \label{lem:perfect-A-stratum}
The map $H_\G^\ast(X_{d_2/2}^\ast\cup \Scal_a, X^\ast_{d_2/2})\to
H_\G^\ast(X_{d_2/2})$ is injective.
\end{lemma}

\begin{proof}
By Proposition \ref{prop:bott1},
 it suffices to show that $H^\ast_\G(\nu^-_a,\nu_a')\to H^\ast_\G(\nu^-_a)$ is injective,
or equivalently, that $H^\ast_\G(\nu^-_a)\to H^\ast_\G(\nu_a')$ is surjective. 
Consider the following commutative diagram:
\begin{equation*}
\xymatrix{
 H^\ast(B\Gcal) \ar[d] \ar[dr]^{\pi^\ast} & \\
H^\ast_\G(\nu_a^-) \ar[r]^{j} & H^\ast_\Gcal(\nu_a') 
}
\end{equation*}
By Lemma \ref{lem:neg-a-cohomology} and  \cite{WentworthWilkin11}, 
$\pi^\ast$ is surjective.  Therefore $j$ is surjective as well.
\end{proof}

Next, we need the following lemma.
\begin{lemma}\label{lem:poincare-polynomials}
Let $(A, B, C)$ be a triple of topological spaces, and suppose that the map $H^*(A, C) \rightarrow H^*(A)$ is injective. Then
$
P_t(A) - P_t(B) = P_t(A, C) - P_t(B, C) 
$.
Moreover, if we suppose in addition that the inclusion of pairs $(B, C) \hookrightarrow (A, C)$ induces a surjection $H^*(A, C) \rightarrow H^*(B, C)$ in cohomology, then the map $H^*(A) \rightarrow H^*(B)$ is a surjection.
\end{lemma}

\begin{remark}
If the inclusions $C \hookrightarrow B \hookrightarrow A$ are inclusions of $G$-spaces, then the above
result is also true in $G$-equivariant cohomology.
\end{remark}

\begin{proof}
We have the following commutative diagram of exact sequences
\begin{equation}\label{eqn:lemma-comm-diagram}
\xymatrix{
\cdots \ar[r] & H^\ast(A, C) \ar[r] \ar[d] & H^\ast(A) \ar[r] \ar[d] & H^\ast(C) \ar[r]
\ar[d]^{\begin{sideways}$\sim$\end{sideways}} & \cdots \\
\cdots \ar[r] & H^\ast(B, C) \ar[r] & H^\ast(B) \ar[r] & H^\ast(C) \ar[r] & \cdots 
}
\end{equation}
The assumption implies that the top horizontal sequence splits, and therefore  the bottom
horizontal sequence also splits. The result follows immediately.
\end{proof}

\begin{proof}[Proof of Theorem \ref{thm:perfection}]
By Lemma \ref{lem:perfect-A-stratum} and Proposition \ref{prop:bott1}, it suffices to consider region $({\bf I})$.
By the argument in \cite[Sect.\ 3.1]{DWWW11}, 
the Atiyah-Bott lemma implies that $H^\ast_\Gcal(X_k, X_k'')\to H^\ast_\Gcal(X_k)$ is injective.
By Corollary \ref{cor:C-bott}, we may then apply Lemma \ref{lem:poincare-polynomials}
 to the triple $(X_k, X_k^\ast, X_k'')$ and conclude that $H^\ast_\Gcal(X_k) \to H^\ast_\Gcal(X_k^\ast)$ surjects.
This completes the proof.  We also record that 
 in this case
\begin{equation} \label{eqn:p1}
P^\Gcal_t(X_\ell)-P^\Gcal_t(X_\ell^\ast)=
P^\Gcal_t(X_\ell, X''_\ell)-P^\Gcal_t(X_\ell^\ast, X''_\ell)
\end{equation}
\end{proof}

\section{The Equivariant Betti Numbers}\label{sec:equivariant-betti-numbers}

\subsection{$\U(2,1)$ bundles} \label{sec:u(2,1)-poincare}  
The calculations in the previous sections lead to the following formula for the equivariant Poincar\'e polynomial of $\mathcal{B}(d_1,d_2)$.
The contributions of individual  strata are as follows.

\begin{enumerate}

\item 
 For the $A$-stratum, use Lemmas \ref{lem:neg-a-cohomology} and  \ref{lem:perfect-A-stratum} to conclude
\begin{align*}
P_t^\Gcal(X^\ast_{d_2/2}\cup\Scal_a)-P_t^\Gcal(X^\ast_{d_2/2}) &=
 \frac{1}{(1-t^2)^2} P_t(J(X)) P_t^\mathcal{G}(\mathcal{A}^{ss}(E_2)) \\
 &\qquad - \frac{1}{(1-t^2)}P_t(\Ncal_{\sigma_{min}}(E_1^* E_2 \otimes K)) P_t(J_{d_1}(X)) \\
& -\begin{cases}
0 & \ \text{if $d_2$ odd} \\
\displaystyle\frac{t^{2(2g-2+d_2/2-d_1)}}{(1-t^2)} P_t(J(X))^2 P_t(S^{2g-2+d_2/2-d_1}X)
& \ \text{if $d_2$ even} 
\end{cases}
\end{align*}

\item   For $
\tfrac{1}{3} (2d_2-d_1) < \ell \leq d_2/2$, \eqref{eqn:mb-c} and \eqref{eqn:c2} imply
\begin{equation*}
P_t^\Gcal(X'_\ell)-P_t^\Gcal(X^\ast_\ell)= \frac{t^{4g-4 + 2 \ell - 2 d_1} }{(1-t^2)^2} P_t(J(X))^2 P_t(S^{\ell - d_1 + 2g-2} X)
\end{equation*}

\item  For $
d_2/2 < \ell \leq \tfrac{1}{3}(d_1+d_2)$, Lemma \ref{lem:b1} and Corollary \ref{cor:b1} imply (recall that $\xi''$ is surjective)
$$
P_t^\Gcal(X_\ell)-P_t^\Gcal(X_\ell^\ast) = \frac{ t^{2(g-1+2\ell - d_2)}}{(1-t^2)^3} P_t(J(X))^3  
  - \frac{t^{2(g-1+2\ell - d_2)} }{(1-t^2)^2} P_t(J(X))^2 P_t(S^{d_2 - d_1 + 2g-2-\ell} X) 
$$

\item For $
 \tfrac{1}{3}(d_1+d_2) < \ell \leq d_2-d_1+2g-2$, it follows from \eqref{eqn:p1}, Lemma \ref{lem:C-bott}, and \eqref{eqn:zeta'_C1}  that
\begin{align*}
P_t^\Gcal(X_\ell)-P_t^\Gcal(X_\ell^\ast)&=
 \frac{t^{2(g-1+2\ell - d_2)}}{(1-t^2)^3}  P_t(J(X))^3  \\
&\qquad - \frac{t^{2(g-1+2\ell - d_2)}}{(1-t^2)}  P_t(J(X)) P_t(S^{d_2-d_1+2g-2-\ell } X) P_t(S^{2g-2-\ell+d_1} X) 
\end{align*}

\item For $
\max\{d_1, d_2-d_1+2g-2\} < \ell \leq d_1+2g-2$, it follows from Proposition \ref{prop:bott1} and eq.'s  \eqref{eqn:c3-bundle-case}, \eqref{eqn:b3}, and \eqref{eqn:nu''_B2}
 that
\begin{align*}
P_t^\Gcal(X'_\ell)-P_t^\Gcal(X_\ell^\ast)&=
 \frac{t^{2(g-1+2\ell - d_2)}}{(1-t^2)^2}  P_t(J(X))^2 P_t(S^{d_1-\ell+2g-2} X)  \\
P_t^\Gcal(X_\ell)-P_t^\Gcal(X'_\ell)&=
 \frac{t^{2(g-1+2\ell - d_2)}}{(1-t^2)^3}  P_t(J(X))^3 
   -\frac{t^{2(g-1+2\ell - d_2)}}{(1-t^2)^2}  P_t(J(X))^2 P_t(S^{d_1-\ell+2g-2} X)
\end{align*}

\item For $
d_1+2g-2 < \ell$, or if $d_2-d_1+2g-2<\ell\leq d_1$,  it follows from Proposition \ref{prop:bott1} and 
\eqref{eqn:b3-large-degree},  from \eqref{eqn:nu'-B1-large} and Remark \ref{rem:b1-c2}, or
 from \eqref{eqn:b2-large-degree},  that
$$
P_t^\Gcal(X_\ell)-P_t^\Gcal(X_\ell^\ast)=
 \frac{t^{2(g-1+2\ell - d_2)}}{(1-t^2)^3}  P_t(J(X))^3 
$$
\end{enumerate}

Applying Theorem \ref{thm:perfection}, we compute
$$
P_t(B \mathcal{G})-P_t^\mathcal{G}(\Bcal^{ss}(d_1, d_2))
=\sum_{\ell\in \Delta_{d_1,d_2}}  P_t^\Gcal(X_\ell)-P_t^\Gcal(X_\ell^\ast) 
$$
Notice that the last term in (i), which occurs only when $d_2$ is even, is exactly canceled by one of the terms in (ii). Combining the remaining terms, we obtain
  
\begin{proposition}\label{prop:u(2,1)-polynomial}
The $\mathcal{G}$-equivariant Poincar\'e polynomial of $\Bcal^{ss}(d_1, d_2)$ is given by
\begin{align} 
P_t^\mathcal{G}(\Bcal^{ss}(d_1, d_2)) &= P_t(B \mathcal{G})
- \frac{1}{(1-t^2)^2} P_t(J(X)) P_t^\mathcal{G}(\mathcal{A}^{ss}(E_2)) -\sum_{d_2/2<\ell}  \frac{t^{2(g-1+2\ell - d_2)}}{(1-t^2)^3}  P_t(J(X))^3\notag \\
& \qquad + \frac{1}{(1-t^2)}P_t(\Ncal_{\sigma_{min}}(E_1^* E_2 \otimes K)) P_t(J_{d_1}(X)) \notag\\
&\qquad +\sum_{d_2/2<\ell\leq \frac{1}{3}(d_1+d_2)}  \frac{t^{2(g-1+2\ell - d_2)} }{(1-t^2)^2} P_t(J(X))^2 P_t(S^{d_2-\ell - d_1 + 2g-2} X) \label{eqn:equiv-poincare-poly}   \\
&\qquad 
-\sum_{\frac{1}{3}(2d_2-d_1) < \ell<d_2/2} \frac{t^{4g-4 + 2 \ell - 2 d_1} }{(1-t^2)^2} P_t(J(X))^2 P_t(S^{\ell - d_1 + 2g-2} X)
\notag \\
&\hskip-1cm +\sum_{\frac{1}{3}(d_1+d_2) < \ell \leq d_2-d_1+2g-2}
  \frac{t^{2(g-1+2\ell - d_2)}}{(1-t^2)}  P_t(J(X)) P_t(S^{2g-2+d_2-\ell - d_1} X) P_t(S^{2g-2-\ell+d_1} X) 
 \notag
\end{align}
\end{proposition}

\begin{proof}[Proof of Theorem \ref{thm:simplified-polynomials}]
We need to show that the expression \eqref{eqn:equiv-poincare-poly} agrees with 
\eqref{eqn:non-coprime-polynomial} and \eqref{eqn:coprime-polynomial}.
By the result of Atiyah-Bott \cite{AtiyahBott83},
$$
P_t(B \mathcal{G})
- \frac{1}{(1-t^2)^2} P_t(J(X)) P_t^\mathcal{G}(\mathcal{A}^{ss}(E_2)) -\sum_{d_2/2<\ell}  \frac{t^{2(g-1+2\ell - d_2)}}{(1-t^2)^3}  P_t(J(X))^3=0
$$
eliminating the first line on the right hand side of  \eqref{eqn:equiv-poincare-poly}.  Let $E=E_1^* E_2 \otimes K$, and recall the definitions \eqref{eqn:sigma} and \eqref{eqn:sigma-min}. 
By 
 \cite[Thm.'s 8.4.1 and 8.4.2]{WentworthWilkin11},
\begin{align*}
P_t^{\Gcal(E)}(\Ccal_{\sigma(d_1,d_2)}(E))
&-
P_t(\Ncal_{\sigma_{min}}(E)) =\\
& +\sum_{d_2/2<\ell< \frac{1}{3}(d_1+d_2)}  \frac{t^{2(g-1+2\ell - d_2)}-t^{2(g-1+d_2 - d_1-\ell)} }{(1-t^2)} P_t(J(X))^2 P_t(S^{d_2-\ell - d_1 + 2g-2} X)   \\
&+\begin{cases} \quad 0 & \text{if $d_1+d_2 \not\equiv 0 \mod 3$} \\
 \displaystyle\frac{t^{2(g-1+\frac{1}{3}(2d_1 - d_2))} }{(1-t^2)} P_t(J(X))^2 P_t(S^{2g-2-\frac{2}{3}(2d_1 - d_2)} X)
 & \text{if $d_1+d_2\equiv 0 \mod 3$}
\end{cases}
\end{align*}
Using this, and 
substituting $\ell\mapsto d_2-\ell$ in the fourth line of  \eqref{eqn:equiv-poincare-poly}, the result follows.
\end{proof}

\begin{proof}[Proof of Corollary \ref{cor:maximal}]
When the Toledo invariant $\tfrac{2}{3}(2d_1 - d_2)$ achieves its maximal value  $2g-2$ then the Poincar\'e polynomial \eqref{eqn:non-coprime-polynomial} simplifies  further. Firstly, note that in this case
$
\tfrac{1}{3}(d_1+d_2)  = d_2-d_1-(2g-2)
$,
and so the summation on the right hand side of \eqref{eqn:non-coprime-polynomial} vanishes.
Secondly, 
for the Bradlow space, $\deg E=g-1=\sigma(d_1,d_2)$. Therefore, in the case of maximal Toledo invariant, the stability parameter  is maximal (and non-generic).  By \cite[Thm.\ 8.4.2]{WentworthWilkin11}, 
the first term on the right hand side of \eqref{eqn:non-coprime-polynomial} is
\begin{align*}
\frac{1}{(1-t^2)} P_t(J(X))^2 P_t(\CBbb P^{2g-3}) & + \frac{t^{4g-4}}{(1-t^2)^2}  P_t(J(X))^2 \\
 & = \frac{1}{(1-t^2)} P_t(J(X))^2 \left( 1 + \cdots + t^{4g-6} + t^{4g-4} (1 + t^2 + \cdots ) \right) \\
 & = \frac{1}{(1-t^2)^2} P_t(J(X))^2 
\end{align*}
\end{proof}

\subsection{$\SU(2,1)$ bundles}

Many of the constructions for $\U(2,1)$ 
Higgs bundles described above also carry over 
to the space $\mathcal{B}_\Lambda(d_1, d_2)$ of 
$\SU(2,1)$ Higgs bundles. In particular, we 
have the same indexing set for the stratification, and the 
index at a critical point can also be computed by an 
analogous calculation to that in Section \ref{sec:neg-directions}.
The major difference between the two cases is that the Kirwan map
$\kappa_0$ from \eqref{eqn:kirwan_map} 
is no longer  necessarily surjective.   However, repeated application of Lemma \ref{lem:poincare-polynomials} allows us to compute the contributions from each critical set individually.

Due to the fixed determinant condition, some of the spaces that contribute to the Poincar\'e polynomial are different to those that appear in the calculation of the previous section: they are finite covers of known spaces (cf. \cite{Hitchin87}, \cite{Gothen94}, \cite{Gothen02}) and so we begin by describing their construction. 

Let $\widetilde{S}(m_1, m_2)$ 
to be the $3^{2g}$-fold cover of $S^{m_1} X \times S^{m_2} X$ 
defined as in the Introduction (see \cite{Gothen94,Gothen02}). These spaces appear in \eqref{eqn:SU-4}.
Recall that the construction is via pullback, as follows
\begin{equation*}
\xymatrix{
\widetilde{S}(m_1, m_2) \ar[r]^{\quad p_2} \ar[d]_{p_1} & J(X) \ar[d]^g \\
S^{m_1} X \times S^{m_2} X  \ar[r]^{\qquad f} & J(X) 
}
\end{equation*}
where $\widetilde{S}(m_1, m_2) \subset S^{m_1} X \times S^{m_2} X \times J(X)$, $p_1$ is projection onto the first two factors, $p_2$ is projection onto the third factor, $f$ maps $(M_1, \varphi_1, M_2, \varphi_2) \mapsto M_1^* M_2 \Lambda$ and $g$ is the three-fold covering map $L \mapsto L^3$. Note that if $M_1 = L_2^* L_3  \otimes K = L_1^* (L_2^*)^2 \Lambda \otimes K$ and $M_2 = L_1^* L_2 \otimes K$, then $M_1^* M_2 \Lambda = L_2^3$, and so $\widetilde{S}(m_1, m_2)$ is the space of bundles $L_1, L_2$ together with nonzero sections $\varphi_1 \in H^0(L_1^* L_2 \otimes K)$ and $\varphi_2 \in H^0(L_1^* (L_2^*)^2 \otimes K)$, where $m_1 = \deg(L_1^* L_2 \otimes K)$ and $m_2 = \deg(L_1^* (L_2^*)^2 \otimes K)$.

As above, let $E=E_1^* E_2 \otimes K$.
For the Type A stratum (see \eqref{eqn:SU-A}), 
define $\widetilde{\Ncal}_\sigma(E)$ 
to be the $3^{2g}$-fold cover of the Bradlow space $\Ncal_\sigma(E)$, which is constructed via the following pullback diagram

\begin{equation*}
\xymatrix{
\widetilde{\Ncal}_\sigma(E) \ar[r]^(0.6){p_2} \ar[d]_{p_1} & J(X) \ar[d]^g \\
\Ncal_\sigma(E) \ar[r]^(0.6)f & J(X) 
}
\end{equation*}
where $p_1(E_1, E_2, \varphi) = (E, \varphi)$, $p_2(E_1, E_2, \varphi) = E_1$, $f(E, \varphi) = \det E$ and $g(L) = (L^*)^3 K^2 \Lambda$. Note that
$$
f \circ p_1(E_1, E_2, \varphi)  = g \circ p_2 (E_1, E_2, \varphi) \
\Longleftrightarrow \  (E_1^*)^2 (\det E_2 ) K^2  = (E_1^*)^3 K^2 \Lambda \
\Longleftrightarrow \ \det(E_2 \oplus E_1)  = \Lambda 
$$

The construction is analogous to \cite[Proposition 2.9]{Gothen94}, but 
the underlying space is different, since 
we use a different stability parameter in this calculation to 
that used in Gothen's calculation ($\sigma_{min}$ as opposed to $\sigma(d_1, d_2)$).
Note, however, that by \cite{Gothen02}, or the methods of \cite{WentworthWilkin11}, it still follows that 
$$
P_t(\widetilde{\Ncal}_{\sigma_{min}}(E))=
P_t({\Ncal}_{\sigma_{min}}(E)) 
$$

The final case to consider is where there are three 
line bundles $L_1, L_2, L_3$ satisfying $L_1 L_2 L_3 = 
\Lambda$, and one section $\varphi \in H^0(L_j^* L_k \otimes K) \setminus 
\{ 0 \}$, where $j, k \in \{ 1, 2, 3\}$ and $j \neq k$. These spaces appear in \eqref{eqn:SU-2}, \eqref{eqn:SU-3} and \eqref{eqn:SU-5} as the cohomology of the type $C$ critical sets and also in \eqref{eqn:SU-A} (when $d_2$ is even).
Let $i \in \{ 1, 2, 3 \} \setminus \{ j, k \}$, and note 
that the fixed determinant condition $L_1 L_2 L_3 = \Lambda$ can be resolved by setting $L_i = \Lambda L_j^* L_k^*$. Then the space under consideration becomes
\begin{align*}
\bigl\{ (L_1, L_2, L_3, \varphi) \, : \, L_1 L_2 L_3 &= \Lambda, \varphi \in H^0(L_j^* L_k \otimes K) \setminus \{ 0 \} \bigr\} \\
 &= \left\{ (L_j, L_k, \varphi) \, : \, \varphi \in H^0(L_j^* L_k \otimes K) \setminus \{ 0 \} \right\} 
\end{align*}
which fibers over $J(X) \times S^{2g-2 + \deg L_k - \deg L_j} X$ with fiber $\CBbb^*$. In particular, if $S^1$ acts freely on the $\CBbb^*$ factor, then the $S^1$-equivariant Poincar\'e polynomial is $P_t(J(X)) P_t(S^{2g-2+\deg L_k - \deg L_j} X)$.

In the same way as for the $\U(2,1)$ case, we can calculate the contributions of the individual strata. These contributions are listed below.

\begin{enumerate}

\item 
 For the $A$-stratum
\begin{align}\label{eqn:SU-A}
\begin{split}
P_t^\Gcal(X^\ast_{d_2/2}\cup\Scal_a)-P_t^\Gcal(X^\ast_{d_2/2}) &=
 \frac{1}{1-t^2} P_t^\mathcal{G}(\mathcal{A}^{ss}(E_2)) 
 - P_t({\Ncal}_{\sigma_{min}}(E_1^* E_2 \otimes K))  \\
& -\begin{cases}
0 & \ \text{if $d_2$ odd} \\
\displaystyle t^{2(2g-2+d_2/2-d_1)} P_t(J(X)) P_t(S^{2g-2+d_2/2-d_1}X)
& \ \text{if $d_2$ even} 
\end{cases}
\end{split}
\end{align}

\item   For $
\tfrac{1}{3} (2d_2-d_1) < \ell \leq d_2/2$
\begin{equation}\label{eqn:SU-2}
P_t^\Gcal(X'_\ell)-P_t^\Gcal(X^\ast_\ell)= \frac{t^{4g-4 + 2 \ell - 2 d_1} }{(1-t^2)^2} P_t(J(X)) P_t(S^{\ell - d_1 + 2g-2} X)
\end{equation}

\item  For $
d_2/2 < \ell \leq \tfrac{1}{3}(d_1+d_2)$ (or $d_2-d_1+2g-2<\ell<d_1$)
\begin{equation}\label{eqn:SU-3}
P_t^\Gcal(X_\ell)-P_t^\Gcal(X_\ell^\ast) = \frac{ t^{2(g-1+2\ell - d_2)}}{(1-t^2)^2} P_t(J(X))^2   - \frac{t^{2(g-1+2\ell - d_2)} }{1-t^2} P_t(J(X)) P_t(S^{d_2-\ell - d_1 + 2g-2} X) 
\end{equation}

\item For $
 \tfrac{1}{3}(d_1+d_2) < \ell \leq d_2-d_1+2g-2$
\begin{align}\label{eqn:SU-4}
\begin{split}
P_t^\Gcal(X_\ell)-P_t^\Gcal(X_\ell^\ast)&=
 \frac{t^{2(g-1+2\ell - d_2)}}{(1-t^2)^2}  P_t(J(X))^2  \\
&\qquad - t^{2(g-1+2\ell - d_2)}   P_t(\widetilde{S}(2g-2+d_2-\ell - d_1, 2g-2-\ell+d_1) )
\end{split}
\end{align}

\item For $
\max\{d_1, d_2-d_1+2g-2\} < \ell \leq d_1+2g-2$
\begin{align}\label{eqn:SU-5}
\begin{split}
P_t^\Gcal(X'_\ell)-P_t^\Gcal(X_\ell^\ast)&=
 \frac{t^{2(g-1+2\ell - d_2)}}{1-t^2}  P_t(J(X)) P_t(S^{d_1-\ell+2g-2} X)  \\
P_t^\Gcal(X_\ell)-P_t^\Gcal(X'_\ell)&=
 \frac{t^{2(g-1+2\ell - d_2)}}{(1-t^2)^2}  P_t(J(X))^2 
  -\frac{t^{2(g-1+2\ell - d_2)}}{1-t^2}  P_t(J(X)) P_t(S^{d_1-\ell+2g-2} X)
\end{split}
\end{align}

\item For $
d_1+2g-2 < \ell$, or if $d_2-d_1+2g-2<\ell\leq d_1$,
\begin{equation}\label{eqn:SU-6}
P_t^\Gcal(X_\ell)-P_t^\Gcal(X_\ell^\ast)=
 \frac{t^{2(g-1+2\ell - d_2)}}{(1-t^2)^2}  P_t(J(X))^2 
 \end{equation}
\end{enumerate}

Theorem   \ref{thm:su-simplified-polynomials} then follows as in the non-fixed determinant case described in the previous section.  We omit the details.

\section{Action of $\Gamma_3$ and the Torelli group} \label{sec:Torelli}

We first 
fix the following   notation. 
Recall that $\Gamma_3=H^1(M, \ZBbb/3)$. 
Then as elements of  $\Gamma_3$ are homomorphisms $\pi\to \ZBbb/3$, $\Gamma_3$
acts on $\Hom(\pi, \SU(2,1))$ by multiplication. 
 The Torelli group $\Tor(M)$ acts on $\Hom(\pi, \SU(2,1))\doubleslash\SU(2,1)$ by  
outer automorphisms of $\pi$.  This induces an action on equivariant
cohomology which commutes with $\Gamma_3$.
 In this section we compute the induced action
of $\Gamma_3\times \Tor(M)$  on the 
$\SU(2,1)$-equivariant cohomology of $\Hom(\pi, \SU(2,1))$.

Following \cite{Looijenga97}, 
   let
$
Q(\Gamma_3)=\{\text{cyclic quotients of $\Gamma_3$}\}
$.
Then $C\in Q(\Gamma_3)$ 
is either $\{0\}$ or $\ZBbb/3$. 
A choice of embedding $ C\hookrightarrow\overline\QBbb$ 
gives a homomorphism $\ZBbb[\Gamma_3]\to \overline\QBbb$, and we 
let $I_C$ denote the kernel.  If
$R_C=\ZBbb[\Gamma_3]/I_C$ and  $K_C=\QBbb\otimes R_C$, then $K_C=\QBbb$ 
if $C=\{0\}$, and otherwise $K_C\cong \QBbb(\xi)$, for $\xi$  a nontrivial third root of unity (though not canonically so).
The field $K_C$ has a natural ``complex conjugation" induced by 
$$
\overline{\sum_{g\in \Gamma_3}c_g g}=\sum_{g\in \Gamma_3}c_g g^{-1}
$$
If $W$ is a $K_C$-vector space, let $\overline W$ denote the vector space with the same underlying $\QBbb$-structure, but where  multiplication by scalars $\lambda\in K_C$ is given by $
\lambda\cdot w=\bar\lambda w$.

Every $\{0\}\neq C\in Q(\Gamma_3)$ gives rise to a connected, unramified 
$3$-fold covering $X_C\to X$.  
Namely, the choice of basepoint $p$ gives an Abel mapping
$X\hookrightarrow J(X)$.  Let $\widetilde X_3$ be the
$\Gamma_3$-covering obtained by pulling back the $\Gamma_3$-covering
$J(X)\to J(X) : L\mapsto L^3$.
Let $X_C$ be
the quotient of $\widetilde X_3$ by the kernel of $\Gamma_3\to C$.
Then $\Gamma_3$ acts on $X_C$ by  
deck transformations,  and there is a decomposition
$$
H^1(X_C,\QBbb)\cong H^1(X,\QBbb)\oplus\left\{ R_C\otimes_{\ZBbb[\Gamma_3]}H^1(X_C,\QBbb) \right\}
$$
where $W_C(X)=R_C\otimes_{\ZBbb[\Gamma_3]}H^1(X_C,\QBbb)$ is a $K_C$-vector space of dimension $2g-2$.
Lifting elements of the  Torelli group  then gives a  surjection of $\Tor(X)$ onto the group of projective unitary transformations of $W_C(X)$, where the unitary structure is the extension by $K_C$ of the symplectic pairing (see \cite{Looijenga97}).

For integers $m_1,m_2\geq 0$, define the $\QBbb$-vector space
\begin{equation} \label{eqn:vm}
V(m_1,m_2)=
\bigoplus_{\{0\}\neq C\in Q(\Gamma_3)} \wedge^{m_1} \overline {W_C(X)}\otimes_{K_C} \wedge^{m_2}{W_C(X)}
\end{equation}
(the exterior products are over $K_C$).
Also,
recall the space $\widetilde S(m_1,m_2)$ from the previous section.
 For $\SU(2,1)$ representations of $\pi$, the 
Toledo invariant is an even integer, and so
$m_1\equiv m_2 \negthinspace \mod 3$. Hence, the diagonal action of $\Gamma_3$ is trivial  on the terms in $V(m_1,m_2)$.
In particular, the projective representation of the Torelli group lifts to a linear one.
With this notation we state

\begin{proposition} \label{prop:decomposition}
The $\Gamma_3$ decomposition is given by 
$$
H^p(\widetilde S(m_1,m_2))
= H^p(S^{m_1}X\times S^{m_2}X) \oplus
\begin{cases} \{0\} & \text{ if $p\neq m_1+m_2$} \\
 V(m_1, m_2) & \text{ if $p=m_1+m_2$}
\end{cases}
$$
\end{proposition}

\begin{proof}
Let $\widetilde S^mX$ be the pull-back of the fibration $S^mX\to J(X)$ under the $3^{2g}$-fold covering $J(X)\to J(X): L\mapsto L^3$.  Then $\Gamma_3$ acts on $\widetilde S^mX$, and by \cite{Looijenga97} we have
$$
H^\ast(\widetilde S^mX, \QBbb)\cong \bigoplus_{C\in Q(\Gamma_3)}  R_C\otimes_{\ZBbb[\Gamma_3]}H^\ast(\widetilde S^mX,\QBbb) 
$$
For $C=\{0\}$ this amounts to
$$
R_C\otimes_{\ZBbb[\Gamma_3]}H^\ast(\widetilde S^mX,\QBbb) =
\left[H^\ast(\widetilde S^mX,\QBbb)\right]^{\Gamma_3}=
H^\ast(S^mX,\QBbb)
$$
For $C\neq \{0\}$, we have an identification of $K_C$-vector spaces
$$
R_C\otimes_{\ZBbb[\Gamma_3]}H^\ast(\widetilde S^mX,\QBbb) \cong
H^\ast(S^mX,{\mathcal F}^{(m)}_C)
$$
where ${\mathcal F}^{(m)}_C\to S^mX$ is a rank-1 local system.  It follows exactly as in Hitchin \cite{Hitchin87} that
there is a rank-1 local system ${\mathcal F}_C\to X$,  such that
$$
H^p(X,{\mathcal F}_C)\cong\begin{cases} \{0\} & \ p=0,2 \\
 W_C(X) & \  p=1
\end{cases}
$$
$$
H^p(S^mX,{\mathcal F}^{(m)}_C)\cong
\begin{cases} \{0\} & \ p\neq m \\
 \wedge^m H^1(X,{\mathcal F}_C) & \ p=m
\end{cases}
$$
Explicitly,  if $\pr:X_C\to X$ is the  covering, the
presheaf  ${\mathcal F}_C(U)$ is given by 
locally constant functions 
$\varphi: \pr^{-1}(U)\to K_C$ 
satisfying
$\varphi(gx)=g\varphi(x)$  for all $x\in X_C$, $g\in \Gamma_3
$.
In the case where the map $S^mX\to J(X)$ is factored through $L\mapsto L^\ast$, then 
$$
H^m(S^mX,{\mathcal F}^{(m)}_C)\cong \wedge^m H^1(X,{\mathcal F}_C^\ast)
$$
and clearly $H^1(X,{\mathcal F}_C^\ast )\cong \overline{W_C(X)}$. 
Applying this argument to $S(m_1, m_2)$, we have
$$
H^\ast(S(m_1,m_2),\QBbb)=
H^\ast(S^{m_1}X\times S^{m_2}X,\QBbb)\oplus
\bigoplus_{\{0\}\neq C\in Q(\Gamma_3)}R_C\otimes_{\ZBbb[\Gamma_3]}H^\ast(S(m_1,m_2),\QBbb)
$$
Now by the Kunneth formula, for $C\neq\{0\}$,
\begin{align*}
R_C\otimes_{\ZBbb[\Gamma_3]}H^p(S(m_1,m_2),\QBbb)
&=
H^p(S^{m_1}X\times S^{m_2}X, {\mathcal F}^{(m_1)}_C\boxtimes
{\mathcal F}^{(m_2)}_C) \\
&=
\bigoplus_{j+k=p}
H^j(S^{m_1}X, {\mathcal F}^{(m_1)}_C)
\otimes_{K_C} 
H^k(S^{m_2}X, {\mathcal F}^{(m_2)}_C)\\ \\
&=
\begin{cases}
\{0\}& \ p\neq m_1+m_2 \\
V(m_1,m_2) &\ p=m_1+m_2
\end{cases}
\end{align*}
\end{proof}

Since  $[K_C: \QBbb]=2$ for $C\neq\{0\}$, and $\# Q(\Gamma_3)=1+\tfrac{1}{2}(3^{2g}-1)$, we have the following
\begin{corollary}[{\cite[Proposition 3.11]{Gothen94}}]  \label{cor:gothen}
If either $m_1$ or $m_2> 2g-2$, then 
$$
P_t(\widetilde S(m_1,m_2))=P_t(S^{m_1}X)P_t(S^{m_2}X)
$$
If $0\leq m_1, m_2\leq 2g-2$,
then
$$
P_t(\widetilde S(m_1,m_2))=P_t(S^{m_1}X)P_t(S^{m_2}X)+(3^{2g}-1){2g-2\choose m_1}{2g-2\choose m_2} t^{m_1+m_2}
$$
\end{corollary}

We now state the result on the action of the Torelli group.
\begin{theorem} \label{thm:torelli}
 Fix a Toledo invariant $0\leq \tau\leq 2g-2$.
Let 
$$
S_\tau=\left\{ 6g-6+\tau/2+2\ell : \ell\in \ZBbb\ ,\ \max\{1,\tau/2\}\leq\ell\leq 2g-2-\tau\right\}$$
Then the  following hold.
\begin{enumerate}
\item The $\SU(2,1)$-equivariant cohomology of $\Hom_\tau(\pi_1(X),\SU(2,1))$
is $\Gamma_3\times\Tor(X)$-invariant in all dimensions $p\not\in S_\tau$.
\item For $p=6g-6+\tau/2+2\ell\in S_\tau$, the nontrivial part of the 
action of $\Gamma_3\times \Tor(X)$
on the $\SU(2,1)$-equivariant cohomology of $\Hom_\tau(\pi_1(X),\SU(2,1))$ in dimension $p$
 is precisely 
$V(m_1,m_2)$, where $m_1=2g-2-\tau-\ell$, $m_2=2g-2+\tau/2-\ell$. 
\end{enumerate}
\end{theorem}

\begin{proof}
Using the stratification $\{X_\ell\}$, 
the argument is the same as in \cite{DWW10}. Note
that  the action on the cohomology of the Bradlow spaces is trivial,
 since Kirwan surjectivity holds for these by \cite{WentworthWilkin11}.
\end{proof}



\end{document}